\documentclass[12pt,reqno]{amsart}
\usepackage[v2,all]{xy}

\usepackage{graphicx} 
\graphicspath{ {FCPic/} }
\usepackage{amsfonts}
\usepackage{mathrsfs}
 \usepackage{amssymb}
 \usepackage{pifont}

 \usepackage[in]{fullpage}

 \usepackage{bold}
 \def\unit{\Eins}
 \def\gh{\mathbbnew{\Gamma}}

 \input{epsf.sty}
 \catcode `\@=11
\def\numberbysection{\@addtoreset{equation}{section}
         \renewcommand{\theequation}{\thesection.\arabic{equation}}}
\numberbysection
\def\subsubsection{\@startsection{subsubsection}{3}%
  \normalparindent{.5\linespacing\@plus.7\linespacing}{-.5em}%
  {\normalfont\bfseries}}

 \newtheorem{thm}{Theorem}[section]
 \newtheorem{lem}[thm]{Lemma}
 \newtheorem{lemma}[thm]{Lemma}
 \newtheorem{prop}[thm]{Proposition}

 \theoremstyle{definition}

 \newtheorem{df}[thm]{Definition}
 \newtheorem{definition}[thm]{Definition}
 
 \newtheorem{rmk}[thm]{Remark}

 \newtheorem{ex}[thm]{Example}

 \newcommand{\leftsub}[2]{{\vphantom{#2}}_{#1}{#2}}

 \def\eps{\epsilon}
 \def\G{\Gamma}
 
 \def\Vect{\mathcal{V}ect}
 \def\dgVect{dg\Vect}

 \def\Fprelie{{\FF}_{pre-Lie}}

 \def\CalC{{\mathcal C}}
 \def\CalD{{\mathcal D}}
 \def\CalE{{\mathcal E}}

 \def\Agg{{\mathcal A}gg}
 \def\Set{{\mathcal S}et}
 \def\Top{{\mathcal T}op}

 \def\Crl{{\mathcal C}rl}

 \def\CCyclic{\mathfrak  {C}}
 
 \def\operads{{\mathfrak O}}
 \def\props{{\mathfrak P}}
 \def\properads{{\mathfrak P}^{ctd}}
 \def\dioperads{{\mathfrak D}}
 \def\modular{{\mathfrak M}}

 \def\F{\mathcal F}
 \def\FF{\mathfrak F}
 
 \def\GG{\mathfrak G}

 \def\FFV{\FF_{\V}}
 
 \def\fL{\mathfrak{L}}

 \def\C{\CalC}
 \def\Z{{\mathbb Z}}
 \def\N{{\mathbb N}}
 \def\G{\Gamma}

 \def\del{\partial}
 \def\colim{\mathrm{colim}}

 \def\Surj{\mathcal{S}urj}
 \def\FSurj{\mathfrak{Sur}}

 \def\O{{\mathcal O}}
 \def\P{{\mathcal P}}
 
 \def\SS{{\mathbb S}}

 \def\odo{\otimes \cdots \otimes}
 \def\K{\mathfrak K}

 \def\day{\circledast}

 \newcommand\ccirc[2]{\, \leftsub{#1}{\circ}_{#2}}
 \def\scirct{\ccirc{s}{t}}
 
 \newcommand\mge[2]{\, \leftsub{#1}{\boxminus}_{#2}}

 \newcommand{\ds}{\displaystyle}

 \def\V{\asts}
 \def\asts{{\mathcal V}}
 \def\F{\clusters}
 \def\clusters{{\mathcal F}}
 
 \def\Ab{{\mathcal Ab}}
 \def\CalG{{\mathcal G}}
 \def\opcat{{\mathcal O }ps}
 \def\op{\mathcal}
 
 \def\smodcat{{\mathcal M}ods}
 \def\forget{\mathit {forget}}
 
 \def\T{{\mathbb T}}
 
 \def\D{\mathcal D}

 \def\oper{op}
 \def\opers{\opcat}

 \def\FT{\mathsf{FT}}
 \def\fops{\F\text{-}\opcat}
 \def\foddops{\F^{odd}\text{-}\opcat}
 \def\fopsc{\F\text{-}\opcat_\C}
 
 \def\vmodsc{\V\text{-}\smodcat_\C}
 \def\vseq{\V\text{-}\text{Seq}_\C}

 \def\Cobar{\Omega}
 \def\Bar{\mathsf{B}}
 \def\tensor{\otimes}
 \def\fopst{\F\text{-}\opcat_{\Top}}

 \newcommand{\id}{\mathbb{I}}

 \newcommand{\fr}{\mathfrak}

 \newcommand{\eops}{\op{E}\text{-}\opcat_\op{C}}

 \newcommand{\vemods}{\V_\fr{E}\text{-}\smodcat_\op{C}}
 \newcommand{\vfmods}{\V_\FF\text{-}\smodcat_\op{C}}

 \newcommand{\adj}[4]{#1\negmedspace: #2\rightleftarrows #3:\negmedspace #4}

 \newcommand{\Fpair}{\F_{dec\O}}
 
 \newcommand{\Fepair}{\Fe_{dec\O}}
 \newcommand{\Fepairtwo}{\Fe'_{dec f_{\ast}(\O)}}

 \newcommand{\Fe}{\mathfrak{F}}

 \def\FFdeco{\FF_{dec\O}}
 \def\Fdeco{\F_{dec\O}}

 \def\Po{\P}
 \def\final{\mathcal T}

\def\B{\mathcal B}

\def\oper{op }
\def\opers{ops }
 \def\Coop{\check\O}
\def\trivial{{\bf 1}}

\def\alephsym{{\mathbb S}}
\def\alephnsym{{\underline{\bf N}}}

\def\Finset{\mathcal{F}inSet}
\def\fopcatc{\F\text{-}\opcat_C}
\def\fopcat{\F\text{-}\opcat}

\def\H{\mathcal H}

\title[Feynman categories]{Lectures on Feynman categories\\[1mm]
{\Small Workshop on higher structures at MATRIX Melbourne\\
 in Creswick, June  7 and 9, 2016}}
\author{Ralph M. Kaufmann}
 \address{Purdue University Department of Mathematics, 150 N.\ University St.,
  West Lafayette, IN 47907}
\email{rkaufman@purdue.edu}

\begin{document}

\begin{abstract}
These are expanded lecture notes from lectures given at the Workshop on higher structures  at MATRIX Melbourne. These notes give an introduction to Feynman categories and their applications.

Feynman categories give a universal categorical way to encode operations and relations. This includes the aspects of operad--like theories such as PROPs, modular operads, twisted (modular) operads, properads, hyperoperads and their colored versions.
There is more depth to the general theory as it applies as well to algebras over operads and an abundance of other related structures, such as crossed simplicial groups, the augmented simplicial category or FI--modules.
Through decorations and transformations the theory is also related to the geometry of moduli spaces. Furthermore the morphisms in a Feynman category give rise to Hopf-- and bi--algebras with examples coming from topology, number theory and quantum field theory. All these aspects are covered.
\end{abstract}

\maketitle

%
\tableofcontents

\section*{Introduction}

\setcounter{section}{1}
\renewcommand{\thesection}{\roman{section}}

\subsection{Main Objective}
Provide a {\it lingua universalis} for operations and relations in order to understand their structure.
The main idea is just like what Galois realized for groups. Namely, one should separate the theoretical structure from the concrete realizations and
representations. What is meant by this is worked out below in the Warm Up section.

In what we are considering, we even take one more step back, namely we provide a theoretical structure for theoretical structures. Concretely the theoretical structures are encoded by a Feynman category and the representations are realized as functors from a given Feynman category $\FF$ to a target category $\C$.
It turns out, however, that to a large extent there are constructions which pass up and down the hierarchy of  theoretical structure vs.\ representation. In concrete examples,
we have a Feynman category whose representations in $\C$ are say algebras. Given a concrete algebra, then there is a new Feynman category whose functors correspond to representations of the algebra. Likewise, for operads, one obtains algebras over the operad as functors.

This illustrates the two basic strategies for acquiring new results. The first is that once we have the definition of a Feynman category,
we can either analyze it further and obtain internal applications to the theory by building several constructions and getting further higher structures. The second is to apply the found results to concrete settings by choosing particular representations.

\subsubsection{Internal Applications} Each of these will be discussed in the indicated section.

\begin{enumerate}
\renewcommand{\theenumi}{\roman{enumi}}
\item Realize universal constructions (e.g.\ free, push--forward, pull--back, plus construction, decorations); see  \S\ref {decopar} and \S\ref{univmasterpar}.
\item Construct universal transforms (e.g.\ bar, co--bar) and model category structure; see \S\ref{modelpar}.
\item Distill universal operations in order to understand their origin (e.g.\ Lie brackets, BV operators, Master equations); see \S\ref{univmasterpar}.
\item Construct secondary objects, (e.g.\ Lie algebras, Hopf algebras); see \S\ref{univmasterpar} and \S\ref{Hopfpar}.
\end{enumerate}
\subsubsection{Applications}
These are mentioned or discussed in the relevant sections and in \S\ref{geopar}.
\begin{enumerate}
\renewcommand{\theenumi}{\roman{enumi}}
 \item Transfer to other areas such as algebraic geometry, algebraic topology, mathematical physics, number theory.
\item Find out information of objects with operations. E.g.\ Gromov-Witten invariants, String Topology, etc.
\item Find out where certain algebra structures come from naturally: pre-Lie, BV, etc.
 \item Find out origin and meaning  of (quantum) master equations.
 \item Construct moduli spaces and compactifications.
 \item Find background for certain types of  Hopf algebras.
 \item Find formulation for TFTs.
 \end{enumerate}

\subsection{References}
The lectures are based on the following references.
\begin{enumerate}
\item with B.~Ward. {\it Feynman categories} \cite{feynman}.
\item with J.~Lucas. {\it Decorated Feynman categories} \cite{deco}.
\item with B.~ Ward. and J.~Zuniga.  {\it The odd origin of Gerstenhaber brackets, Batalin--Vilkovisky operators  and master equations} \cite{KWZ}.
\item with I. Galvez--Carrillo and A. Tonks. {\it Three Hopf algebras and their operadic and categorical background} \cite{GKT}.
\item with C.~Berger. {\it Derived Feynman categories and modular geometry} \cite{BK}.
\end{enumerate}
We also give some brief information on works in progress \cite{KKreimer} and further developments \cite{Ward17}.

\subsection{Organization of the Notes}
These notes are organized as follows. We start with a warm up in \S\ref{warmpar}.
This explains how to understand the concepts mentioned in the introduction. That is, how to construct the theoretical structures in the basic examples of group representations and associative algebras. The section also contains a glossary of the terms used in the following. This makes the text more self--contained. We give the most important details here, but refrain from the lengthy full fledged definitions, which can be found in the standard sources.

In \S\ref{fdefpar}, we then give the definition of a Feynman category and provide the main structure theorems, such as the monadicity theorem and the theorem establishing push--forward and pull--back. We then further explain the concepts by expanding the notions and providing details. This is followed by a sequence of examples. We also give a preview of the examples of operad--like structures that are discussed in detail in \S\ref{zoopar}. We end \S\ref{fdefpar}  with a discussion of the connection to physics and a preview of the various constructions for Feynman categories studied in later sections.

 \S\ref{zoopar} starts by introducing the  category of graphs of Borisov--Manin and the Feynman category $\GG$ which is a subcategory of it. We provide an analysis of this category, which is pertinent to the following sections as a blue print for generalizations and constructions. The usual zoo of operad--like structures is obtained from $\GG$ by decorations and restrictions, as we explain. We also connect the language of Feynman categories to that of operads and operad--like structures. This is done in great detail for the readers familiar with these concepts. We end with omnibus theorems for these structures, which allow us to provide all the three usual ways of introducing these structures (a) via composition along graphs, (b) as algebras over a triple and (c) by generators and relations.

Decoration is actually a technical term, which is explained in \S\ref{decopar}. This paragraph also contains a discussion of so--called non--Sigma, aka.\ planar versions.
We also give the details on how to define the decorations of \S\ref{zoopar} as decorations in the technical sense. We then discuss how with decorations one can obtain the three formal geometries of Kontsevich and end the section with an outlook of further applications of this theory.

The details of enrichments are studied in \S\ref{enrichedpar}. We start by motivating these concepts through the concrete consideration of algebras over operads. After this prelude, we delve into the somewhat involved definitions and constructions. The central ones are Feynman categories indexed enriched over another Feynman category, the $+$ and $hyp$ constructions.
These are tied together in the fact that enrichments indexed over $\FF$ are equivalent to strict symmetric monoidal functors with source $\FF^{hyp}$. This is the full generalization of the construction of the Feynman category for algebras over a given operad. Further constructions are the free monoidal construction $\FF^{\boxtimes}$ for which strict symmetric monoidal functors from
$\FF^{\boxtimes}$ to $\C$ are equivalent to ordinary functors from $\FF$. And the nc--construction $\FF^{nc}$ for which the strict symmetric monoidal functors from $\FF^{nc}$ to $\C$ are equivalent to lax monoidal functors from $\FF$.

Universal operations, transformations and Master equations are treated in \S\ref{univmasterpar}. Examples of universal operations are the pre--Lie bracket for operads or the BV structure for non--connected modular operads. These are also the operations that appear in Master Equations. We explain that these Master Equations are  equations which appear in the consideration of Feynman transforms. These are similar to bar-- and cobar constructions that are treated as well. We explain that the fact that the universal operations appear in the Master Equation is not a coincidence, but rather is a reflection of the construction of the transforms. The definition of the transforms involves odd versions for the Feynman categories, the construction of which is also spelled out.

As for algebras, the bar--cobar or the double Feynman transformation are expected to give resolutions. In order to make these statements precise, one needs a Quillen model structure. These model structures are discussed in \S\ref{modelpar} and we give the conditions that need to be satisfied in order for the transformations above to yield a cofibrant replacement. These model structure are on categories of strict symmetric monoidal functors from the Feynman category into a target category $\C$. The conditions for $\C$ are met for simplicial sets, dg--vector spaces in characteristic $0$ and for topological spaces. The latter requires a little extra work. We also give a W--construction for the topological examples.

 The geometric counterpart to some of the algebraic constructions is contained in \S\ref{geopar}. Here we show how the examples relate to various versions of moduli spaces and how master equations correspond to compactifications.

 Finally, in \S\ref{Hopfpar} we expound the connection of Feynman categories to Hopf algebras. Surprisingly, the examples considered in \S\ref{fdefpar} already yield Hopf algebras that are fundamental to number theory, topology and physics. These are the Hopf algebras of Goncharov, Baues and Connes--Kreimer. We give further generalizations and review the full theory.

\subsection{Acknowledgements}
I thankfully acknowledge my co--authors with whom it has been a pleasure to work. I furthermore thank the organizers of the MATRIX workshop for providing the opportunity to give these lectures and for arranging the special issue.

The work presented here has at various stages been supported by the Humboldt Foundation, the Institute for Advanced Study, the Max--Planck Institute for Mathematics, the IHES and by the NSF. Current funding is provided by the Simons foundation.

\setcounter{section}{0}
\renewcommand{\thesection}{\arabic{section}}
\section{Warm up and Glossary}
\label{warmpar}
Here we will discuss how to think about operations and relations in terms of theoretical structures and their representations by looking at two examples.

\subsection{Warm up I:  Categorical formulation for representations of a group $G$.}
Let  $\underline G$ the category with one object  $*$ and morphism set $G$.
The composition of morphisms is given by group multiplication $f\circ g:= fg$.
This is associative and has the group identity  $e$ as a unit $e=id_*$.

There is more structure though. Since $G$ is a group, we have the extra structure of inverses.
That is every morphism in $\underline G$ is invertible and hence $\underline G$ is a groupoid.
Recall that a category in which every morphism is invertible is called a groupoid.

\subsubsection{Representations as functors}
A representation $(\rho,V)$ of the group $G$  is equivalent to a functor $\underline\rho$ from $\underline G$ to the category of $k$-vector spaces $\Vect_k$.
Giving the values of the functor on the sole object and the morphisms provides:
$\underline\rho(*)=V$, $\rho(g):=\underline\rho(g)\in Aut(V)$. Functoriality then says $\bar\rho(G)\subset Aut(V)$ is a subgroup and all the
relations for a group representation hold.

\subsubsection{Categorical formulation of Induction and Restriction}
Given a morphism $f:H\to G$ between two groups. There are the restriction and induction of any representation $\rho$: $Res^G_H \rho$ and $Ind_H^G \rho$.
The morphism $f$ induces a functor $\underline f$ from $\underline H$ to $\underline G$ which sends the unique object to the unique object and a morphism $g$ to $f(g)$. In terms of  functors restriction  simply becomes pull--back $\underline f^*(\underline \rho):=\rho\circ \underline f$ while induction becomes push--forward,
$\underline{f}_*$, for functors. These even form an adjoint pair.

%

\subsection
{Warm up II: Operations and Relations. Description of Associative Algebras}
An associative algebra in a tensor category $(\C,\otimes)$ is usually given by the following data:
 An object $A$ and one operation: a multiplication $\mu: A\otimes A\to A$
which satisfies the axiom of the associativity equation:
$$(ab)c=a(bc)$$

\subsubsection{Encoding}

Think of $\mu$ as a 2-linear map. Let $\circ_1$ and $\circ_2$ be substitution in the 1st respectively the 2nd variable. This allows us to
rewrite the associativity equation as $$(\mu\circ_1\mu)(a,b,c):=\mu(\mu(a,b),c)=(ab)c=a(bc)=\mu(a,\mu(b,c)):= (\mu\circ_2\mu)(a,b,c)$$

The associativity hence becomes
\begin{equation}
\label{asseq}
\mu\circ_1\mu=\mu\circ_2\mu
\end{equation}
as morphisms $A\otimes A\otimes A \to A$. The advantage of \eqref{asseq} is that it is independent of elements and of $\C$ and merely uses the fact that in multi--linear functions one can substitute. This allows the realization that  associativity is an equation about iteration.

In order to formalize this, we have to allow all possible iterations. The realization this description affords is that all iterations of $\mu$ resulting in an $n$--linear map
are equal. On elements one usually writes $a_1\odo a_n\to a_1\dots a_n$.

In short: for an associative algebra one has one basic operation and the relation is that all $n$--fold iterates agree.

\subsubsection{Variations} If $\C$ is symmetric, one can also consider the permutation action.
Using elements the permutation action gives the opposite multiplication $\tau\mu(a,b)=\mu\circ \tau (a,b)=ba$.

 This give a permutation action on the iterates of $\mu$. It is a free action  and there are $n!$ $n$--linear morphisms generated by $\mu$ and the transposition.
One can also think of commutative algebras or unital versions.
\subsubsection{Categories and functors}
In order to construct the data, we need to have the object $A$, its tensor powers and the multiplication map.
Let $\trivial$ be the category with one object $*$ and one morphism $id_*$. We have already seen that  the functors from $\trivial$ correspond to objects of $\C$.
To get the tensor powers, we let $\alephnsym$ be the category whose objects are the natural numbers including $0$ with only identity morphisms. This becomes a monoidal category with the tensor product given by addition $m\otimes n=m+n$.
Strict monoidal functors $\O$ from $\alephnsym\to\C$ are determined by their value on $1$. Say $\O(1)=A$ then $\O(n)=A^{\otimes n}$.

To model associative algebras, we need a morphisms $\pi:2\to 1$. A monoidal functor $\O$ will assign a morphism $\mu:=\O(\pi):A\otimes A \to A$.
If we look for the ``smallest monoidal category'' that has the same objects as $\alephnsym$ and contains $\pi$ as a morphism, then this is the category $sk(\Surj_<)$
of order preserving surjections between the sets $\underline{n}$ in their natural order. Here we think of $n$ as $\underline n=\{1,\dots, n\}$. Indeed any such surjection is an iteration of $\pi$. Alternatively, $sk(\Surj_<)$ can be constructed from $\alephnsym$ by adjoining the morphism $\pi$ to the strict monoidal category and modding out by the equation analogous to \eqref{asseq}$:\pi\circ id\otimes \pi=\pi\circ \pi\otimes id$.

It is easy to check that functors from $sk(\Surj_<)$ to $\C$ correspond to associative algebras (aka.\ monoids) in $\C$. From this we already gained that starting from say $k$-algebras, i.e.\ $\C=Vect_k$ (the category of $k$ vector spaces), we
can go to any other monoidal category $\C$ and have algebra objects there.

\subsubsection{Variations}
The variation in which we consider the permutation operations is very important. In the first step, we will need to consider $\alephsym$, which has the same objects as $\alephnsym$, but has additional isomorphisms. Namely $Hom(n,n)=\SS_n$ the symmetric group on $n$ letters. The functors out of $\alephnsym$ one  considers are strict symmetric monoidal functors $\O$ into symmetric monoidal categories $\C$. Again, these are fixed by $\O(1)=:A$, but now every $\O(n)=A^{\otimes n}$ has the $\SS_n$ action of permuting the tensor factors according to the commutativity constraints in $\C$.

Adding the morphisms $\pi$ to $\alephsym$ and modding out by the commutativity equations, leaves the ``smallest symmetric monoidal category" that contains the necessary structure. This is the category of all surjections $sk(\Surj)$ on the sets $\underline{n}$. Functors from this category are commutative algebra objects, since $\pi\circ \tau=\pi$ if
$\tau$ is the transposition.

In order to both have symmetry and not force commutativity, one formally does not mod out by the commutativity equations. The result is then equivalent to the category $sk(\Surj_{ord})$ of ordered finite sets with surjections restricted to the sets $\underline{n}$.
 The objects of $\Surj_{ord}$ are a finite set $S$ with an order $<$. The bijections of $S$ with itself act simply transitively on the orders by push--forward.

The second variation is to add an identity. An identity in a $k$--algebra $A$ is described by an element $1_A$, that is a morphism $\eta: k\to A$ with $\eta(1_k)=1_A$. Coding this means that we will have to
have one more morphism in the source category. Since $k=\unit$ is the unit of the monoidal structure of $Vect_k$, we see that we need a morphism $u:0\to1$. We then need to mod out by the appropriate equations, which are given by $\eta\circ_1\mu=\eta\circ_2\mu=id$ which translate to $\pi\circ u\otimes id_1=\pi \circ id_1\otimes u=id_1$.

\subsection{Observations}
\begin{figure}
    \centering
    \includegraphics[width=.7\textwidth]{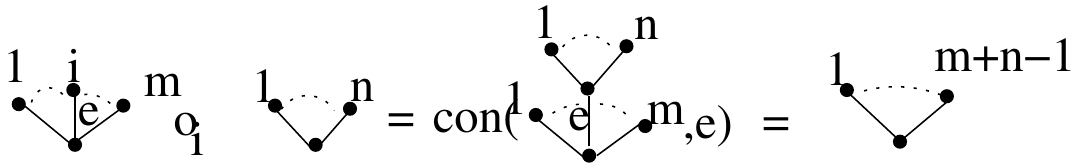}
\caption{\label{corollagraft}
Example of grafting two (planar)  corollas. First graft at a leaf and then contract the  edge.}
\end{figure}

There is a graphical calculus that goes along with the example above. This is summarized in Figure \ref{corollagraft}. Adding in the orders corresponds to regarding planar corollas.

We have dealt with strict structures and actually skeletal structures in the examples. This is not preferable for a general theory. Just as it is preferable to work with all finite dimensional vector spaces in lieu of just considering the collection of $k^n$ with matrices as morphisms.

\subsection{Glossary: Key concepts and Notations} Here is a brief description of key concepts. For more information and full definitions see e.g.\ \cite{MacLane,JoyalStreet}.

{\sc Groupoid:} A category in which every morphism is an isomorphism.

As we have seen, every group defines a groupoid. Furthermore for any category $\C$, the subcategory $Iso(\C)$ which has the same objects as $\C$ but only includes the isomorphisms of $\C$ is a groupoid.

{\sc Monoidal category:} A category $\C$ with a functor $\otimes: \C\times \C \to \C$, associativity constraints and unit constraints. That is an operation on objects $(X,Y)\to X\otimes Y$ and on morphisms $(\phi:X\to Y,\psi:X'\to Y') \to \phi\otimes \psi:X\otimes X'\to Y\otimes Y' $. Furthermore a unit object $\unit$ with isomorphisms $\unit\otimes X\simeq X\simeq X\otimes \unit$ called left and right unit constraints and associativity constraints, which are isomorphisms $a_{X,Y,Z}:X\otimes (Y\otimes Z)\to (X\otimes Y)\otimes Z$. These have to satisfy extra conditions called the pentagon axiom and the triangle equation ensuring the compatibilities. In particular, it is the content of Mac Lane's coherence Theorem that due to these axioms any two ways to iteratively rebracket and add/absorb identities to go from one expression to another are equal as morphisms.

A  monoidal category is called {\em strict}  if the associativity and unit constraints are identities. Again, due to Mac Lane, every monoidal category is   monoidally equivalent to a strict monoidal category (see below).

An example is $\Vect_k$ the category of  $k$-vector spaces with tensor product $\otimes$. Strictly speaking, the associativity constraint $a_{U,V,W}$ acts on
elements as $a_{U,V,W}((u\otimes v)\otimes w))=u\otimes (v\otimes w)$. The unit is $k$ and the unit constraints are $k\otimes U\simeq U\simeq U\otimes k$.

{\sc  Monoidal functor:} A (lax) monoidal functor between two monoidal categories $\C$ and $\D$ is an ordinary functor $F:\C\to \D$ together with a morphisms $\phi_0:\Eins_\D\to F(\Eins_\C)$ and a family of natural morphisms $\phi_2:F(X)\otimes_\D F(Y)\to F(X\otimes_\C Y)$, which satisfy compatibility with associativity and the unit.
A monoidal functor is called strict if these morphisms are identities and strong if the morphism are isomorphisms. If the morphisms go the other way around, the functor is called co-monoidal.

{\sc Symmetric Monoidal Category:}  A  monoidal category $\C$ with all the structures above together with commutativity constraints which are isomorphisms $c_{X,Y}:X\otimes Y \to Y\otimes X$.
These have to satisfy the axioms of the symmetric group, i.e.\  $c_{Y,X}\circ c_{X,Y}=id$ and the braiding for three objects. Furthermore, they are compatible with the associativity constraints, which is expressed by the so--called hexagon equation.

For $\Vect_k$, the symmetric structure $c_{U,V}$ is given on elements as $c_{U,V}(u\otimes v)=v\otimes u$.
We can also consider $\Z$--graded vector spaces. In this category, the commutativity constraint on elements is given by $c_{U,V}(u\otimes v)=(-1)^{deg(u)deg(v)}v\otimes u$ where $deg(u)$ is the $\Z$--degree of $u$.

{\sc Symmetric monoidal functors:} A symmetric monoidal functor is a monoidal functor, for which the $\phi_2$ commute with the  commutativity constraint.

{\sc Free monoidal categories:} There are several versions of these depending on whether one is using strict or non--strict  and symmetric versions or non--symmetric versions.

Let $\V$ be a category. A free (strict/symmetric) category  on $\V$ is a (strict/symmetric) monoidal category $\V^{\otimes}$ and a functor $\jmath:\V\to \V^{\otimes}$ such that
any functor $\imath:\V\to \F$ to a (strict/symmetric) category $\F$ factors as

\begin{equation}
\label{itensordefeq}
\xymatrix{
\V\ar[dr]^{\jmath}\ar^{\imath}[rr]&&\F\\
&\V^{\otimes}\ar^{\imath^{\otimes}}[ur]&
}
\end{equation}
where $\imath^{\otimes}$ is a (strict/symmetric) monoidal functor.

The free strict monoidal category is given by words in objects of $\V$ and words of morphisms in $\V$. The free monidal category is harder to describe. Its objects are iteratively build up from $\otimes$ and the constraints, see \cite{JoyalStreet}, where it is also shown that:
\begin{prop}
 There is a strict monoidal equivalence between the free monoidal category and the strict free monoidal category.
 \end{prop}
This allows us some flexibility when we are interested in data given by a category up to equivalence.

If one includes ``symmetric'' into the  free monoidal category, then one (iteratively) adds morphisms to the free categories that are given by the commutativity constraints. In the strict case, one gets commutative words, but extra morphisms from the commutativity constraints.
As an example,  regard the trivial category $\trivial$: $\trivial^{\otimes,strict}=\alephnsym$ while $\trivial^{\otimes symmetric,strict}=\alephsym$.

{\sc Skeleton of a category:} A skeleton $sk(\C)$ of a category $\C$ is a category that is equivalent to $\C$, but only has one object in each isomorphism class.

An example is the category of ordered finite sets $FinSet$ and morphisms between them with the disjoint union as a symmetric monoidal category.
A skeleton for this category is given by the category whose objects are natural numbers, where each such object $n$ is thought of as the set $\underline n=\{1,\dots,n\}$ and all morphisms between them. This category is known as the (augmented) crossed simplicial group $\Delta_+S$.

{\sc Underlying discrete category:} The underlying discrete category of a category $\C$ is the subcategory which has the same objects as $\C$, but retains the identity maps. It will be denoted by $\C_0$. For instance $\alephsym_0=\alephnsym$.

{\sc Underlying groupoid of a category:} For a category $\C$ the underlying groupoid $Iso(\C)$ is the subcategory of $\C$ which has the same objects as $\C$ buy only retains all the isomorphisms in $\C$.

{\sc Comma categories:}
Recall that for two functors $\imath:\D\to \C$ and $\jmath:\CalE\to \C$, the comma category
$(\jmath\downarrow \imath)$ is the category whose objects are morphisms $\phi\in Hom_{\C}(\jmath(X), \imath(Y))$.
A morphism between such $\phi$ and $\psi$ is given by a commutative diagram.
$$
\xymatrix{
\jmath(X)\ar[d]_{\jmath(f)}\ar[r]^{\phi}&\imath(Y)\ar[d]^{\imath(g)}\\
\jmath(X')\ar[r]_{\psi}&\imath(Y')
}
$$
with $f\in Hom_{\D}(X,X'), g\in Hom_{\CalE}(Y,Y')$. We will write $(\imath(f),\imath(g))$ for such morphisms or simply $(f,g)$.

If a functor, say $\imath: \V\to \F$, is fixed we will just write  $(\F\downarrow \V)$,
and given a category $\CalG$ and an object $X$ of $\CalG$, we denote  the respective comma category by $(\CalG\downarrow X)$. I.e.\ objects are morphisms $\phi: Y\to X$ with $Y$ in $\CalG$ and
morphisms are morphisms over $X$, that is morphisms $Y\to Y'$ in $\CalG$ which commute with the base
maps to $X$. This is sometimes also called the slice category or the category of objects over $X$.

\section{Feynman categories}
\label{fdefpar}
With the examples and definitions of the warm up in mind, we give the definition of Feynman categories and then discuss several basic examples.
The Feynman categories will give the operations and relations part. The concrete examples of the structures thus encoded are then given via functors, just like discussed above.

\subsection{Definition}
\subsubsection{Data for a Feynman category}
\begin{enumerate}
\item $\V$ a groupoid
\item $\F$ a  symmetric monoidal category
\item $\imath: \V\to \F$ a functor.
\end{enumerate}

Let  $\V^{\otimes}$ be the free symmetric category on $\V$ and $\imath^\otimes$ the functor in \eqref{itensordefeq}.
%

\subsection{Feynman category}
\begin{definition}
\label{feynmandef}
The data of triple $\FF=(\V,\F,\imath)$ as above  is called a Feynman category
if the following conditions hold.

\begin{enumerate}
\renewcommand{\theenumi}{\roman{enumi}}

\item \label{objectcond}
 $\imath^{\otimes}$ induces an equivalence of symmetric monoidal categories between $\V^{\otimes}$ and $Iso(\F)$.

\item \label{morcond} $\imath$ and $\imath^{\otimes}$ induce an equivalence of symmetric monoidal categories between $(Iso(\F\downarrow \V))^{\otimes}$ and  $Iso(\F\downarrow \F)$ .

\item For any $\ast\in \asts$, $(\clusters\downarrow\ast)$ is
essentially small.

\end{enumerate}

\end{definition}
Condition (i) is called the {\em isomorphisms condition}, (ii) is  called the {\em hereditary condition} and (iii) the {\em size condition}.
The objects of $(\F\downarrow \V)$ are called {\em one--comma generators}.

\subsubsection{Non-symmetric version}
Now let $(\V,\F,\imath)$ be as above with the exception that $\F$ is only a monoidal category, $\V^{\otimes}$ the free monoidal category, and $\imath^{\otimes}$ is the corresponding morphism of monoidal groupoids.

\begin{df}
A non-symmetric triple $\FF=(\V,\F,\imath)$ as above is called a \emph{non-$\Sigma$} Feynman category  if

\begin{enumerate}
\renewcommand{\theenumi}{\roman{enumi}}

\item \label{nonsymobjectcond}
 $\imath^{\otimes}$ induces an equivalence of  monoidal groupoids between $\V^{\otimes}$ and $Iso(\F)$.

\item \label{nonsymmorcond} $\imath$ and $\imath^{\otimes}$ induce an equivalence of  monoidal groupoids $Iso(\F\downarrow \V)^{\otimes}$ and  $Iso(\F\downarrow \F)$.

\item For any object $\ast_v$ in $\asts$, $(\clusters\downarrow\ast_v)$ is
essentially small.

\end{enumerate}

\end{df}

\subsection{$\opcat$ and $\smodcat$} 

\begin{definition}
Fix   a symmetric monoidal category $\CalC$ and $\FF=(\V,\F,\imath)$  a Feynman category.
\begin{itemize}
\item
$\fopsc:=Fun_{\otimes}(\clusters,\CalC)$ is defined to be the category of strong symmetric monoidal functors
which we will call $\F$--ops in $\CalC$. An object of the category will be referred to as an $\F$-$\oper$ in $\C$.
\item $\vmodsc:=Fun(\asts,\CalC)$, the set of (ordinary) functors will be
called $\V$-mods in $\CalC$ with elements being called a $\V$--mod in $\C$.
\end{itemize}
\end{definition}
There is an obvious forgetful functor $G:\opcat\to \smodcat$ given by restriction.
%

\begin{thm}
The  forgetful functor $G:\opcat\to \smodcat$ has a left adjoint $F$ (free functor) and this adjunction is monadic. This means that the category of
the algebras over the triple $\T=GF$ in $\C$ are equivalent to the category of $\F$-$\opcat_\C$.
\end{thm}

Morphisms between Feynman categories are given by strong monoidal functors that preserve the structures. Natural transformations between them
give 2--morphisms.
The categories $\F$-$\opcat_\C$ and $\F$-$\smodcat_\C$ again are symmetric monoidal categories, where the  symmetric monoidal structure is inherited from $\C$. E.g.\ the tensor product is pointwise, $(\O\otimes \O')(X):=\O(X)\otimes \O'(X)$, and the unit is the functor $\Eins_{\opcat}:\F\to \C$. I.e.\  the functor that assigns $\unit_\C\in Obj(\C)$ to any object in $\V$, and which sends morphisms to the identity morphism.
This is a strong monoidal functor by using the unit constraints.

\begin{thm}
Feynman categories form a 2--category and it has push--forwards and pull--backs for $\opcat$. That is , for a morphism of Feynman categories $f$, both push--forward $f_*$ and pull-back  $f^*$ are adjoint symmetric monoidal functors $f_*:\F-\opcat_C \leftrightarrows \F'-\opcat_C:f^*$.
\end{thm}

\subsection{Details}
\subsubsection{Details on the Definition}
\label{detaildefpar}
The conditions can be expanded and explained as follows.

\begin{enumerate}
\item Since $\V$ is a groupoid, so is $\V^{\otimes}$. Condition (i) on the object level says, that any object $X$ of $\F$ is isomorphic to a tensor product of objects coming from $\V$. $X\simeq \bigotimes_{v\in I}\imath(\ast_v)$. On the morphisms level it  says that all the isomorphisms in $\F$ basically come from $\V$ via tensoring basic isomorphisms of $\V$, the commutativity and the associativity constraints.
 In particular, any two decompositions of $X$ into  $\bigotimes_{v\in I}\imath(\ast_v)$ and $\bigotimes_{v'\in I}\imath(\ast_v')$ there is a bijection $\Psi:I\leftrightarrow I'$ and an isomorphism $\sigma_v:\imath(\ast_v)\to \imath(\ast_{v'})$.
 This implies that for any $X$ there is a unique length $|I|$, where $I$ is any index set for a decomposition of $X$ as above, which we denote by $|X|$. The monoidal unit $\Eins_F$ has length $0$ as the  tensor product over the empty index set.

\item Condition (ii) of the definition of a Feynman category is to be understood as follows: An object in $(\F\downarrow \V)$ is a morphism $\phi:X\to \imath (\ast)$, with $\ast$ in  $Obj(\V)$. An object in $(\F\downarrow \V)^{\otimes}$ is then a formal tensor product of such morphisms, say  $\phi_v: X_v\to \imath(\ast_v)$, $v\in I$ for some index set $I$. To such a formal tensor product, the induced functor assigns $\bigotimes_{v\in V}\phi_v:\bigotimes_v X_v\to \bigotimes_v \ast_v$, which is a morphisms in $\F$ and hence an object of $(\F\downarrow \F)$.

The functor is defined in the same fashion on morphisms. Recall that an isomorphism in a comma category is given by a commutative diagram, in which the vertical arrows are isomorphisms, the horizontal arrows being source and target.
In our case the equivalence of the categories on the object level says that any morphisms $\phi:X\to X'$ in $\F$  has a ``commutative decomposition diagram'' as follows

\begin{equation}
\label{morphdecompeq}
\xymatrix
{
X \ar[rr]^{\phi}\ar[d]_{\simeq}&& X'\ar[d]^{\simeq} \\
 \bigotimes_{v\in I} X_v\ar[rr]^{\bigotimes_{v\in I}\phi_{v}}&&\bigotimes_{v\in I} \imath(\ast_v)
}
\end{equation}
which means that when
$\phi: X\to X'$ and  $X'\simeq \bigotimes_{v\in I} \imath(\ast_v)$ are fixed  there are $X_v\in \F,$ and $\phi_v\in Hom(X_v,\ast_v)$ s.t. the above diagram commutes.

The morphisms part of the equivalence of categories means the following:
\begin{enumerate}
\item For any two such decompositions $\bigotimes_{v\in I} \phi_v$ and $\bigotimes_{v'\in I'}\phi'_{v'}$
there is a bijection
$\psi:I\to I'$ and isomorphisms
$\sigma_v:X_v\to X'_{\psi(v)}$ s.t. $P^{-1}_{\psi}\circ \bigotimes_v \sigma_v\circ \phi_v =\bigotimes \phi'_{v'}$ where
$P_{\psi}$ is the permutation corresponding to $\psi$.
\item These are the only isomorphisms between morphisms.
\end{enumerate}

As it is possible that $X_v=\Eins$, the axiom allows to have morphisms $\Eins\to X'$,
which are decomposable as a tensor product of morphisms $\Eins\to \imath(\ast_v)$. On the other hand, there can be no morphisms $X\to \Eins$ for
any object $X$ with $|X|\geq 1$. If $\Eins$ is the target,  the index set $I$ is empty and hence $X\simeq \Eins$, since the  tensor product over the empty set is the monoidal unit.

We set the length of a morphisms to be $|\phi|=|X|-|X'|$. This can be positive or negative in general. In many interesting examples, it is, however, either non--positive or non--negative.

\item The last condition is a size condition, which ensures that certain colimits over these comma--categories to cocomplete categories exist.
\end{enumerate}

\subsubsection{Details on the adjoint free functor}

The free functor $F$ is defined as follows:
Given a $\asts$--module $\Phi$, we extend $\Phi$
to all objects of $\clusters$ by picking a functor $\jmath$ which yields the equivalence of
$\V^{\otimes}$ and $Iso(\F)$.  Then, if $\jmath(X)=\bigotimes_{v\in I}\ast_v$, we set
\begin{equation}
\Phi(X):=\bigotimes_{v\in I}\Phi(\ast_v)
\end{equation}

Now, for any $X\in\clusters$ we set
\begin{equation}
\label{kandefeq}
F(\Phi)(X)=\colim_{Iso(\clusters\downarrow X)}\Phi\circ s
\end{equation}
where $s$ is the source map in $\clusters$ from $Hom_{\clusters}\to Obj_{\clusters}$ and on the right hand
side or \eqref{kandefeq}, we mean the underlying object. These colimits exist due to condition (iii).
For a  given morphism $X\to Y$ in $\F$, we get an induced morphism of the colimits  and it is straightforward that this defines a functor.
This is actually nothing but the left Kan extension along the functor $\imath^\otimes$ due to (i).
What remains to be proven is that this functor is actually a strong symmetric monoidal functor, that is that  $f_*(\O):\F'\to \C$ is strong symmetric monoidal. This can be shown by using the hereditary condition (ii).

The fact that $f_*$ is itself symmetric monoidal amounts to a direct check as does the fact that $f^*$ and $f_*$ are adjoint functors. The fact that $f^*$ is symmetric monoidal is clear.

\subsubsection{Details on Monadicity}
A triple aka.\ monad on a category is the categorification of a unital semigroup. I.e.\
a triple $\T$ on a category $\C$ is an endofunctor $T:\C\to \C$ together with two natural transformations, $\eta:Id_\C\to T$, where $id_\C$ is the identity functor and a multiplication natural transformation $\mu:T\circ T\to T$, which satisfy the associativity equation $\mu\circ T\mu=\mu\circ \mu T$ as natural transformations  $T^{3}\to T$,
and the unit equation $\mu\circ T\eta=\mu\circ \eta T=id_T$, where $id_T$ is the identity natural transformation of the functor $T$ to itself.
The notation is to be read as follows: $\mu\circ T\mu$ has the components $T(T^2 (X))\stackrel{T(\mu_X)}{\to}T^2X\stackrel{\mu_X}{\to}TX$, where $\mu_X:T^2X\to TX$ is the component of $\mu$.

An algebra over such a triple is an object $X$ of $\C$ and a morphism $h:TX\to X$ which satisfies the unital algebra equations.
$h\circ Th=h\circ \mu_X:T^2X\to X$ and $id_X=h\circ \eta_X:X\to X$.

\subsubsection{Details on morphisms, push--forward and pull--back}
A morphisms of Feynman categories $(\V,\F,\imath)$ and $(\V',\F',\imath')$ is a pair of functors $(v,f)$  where $v\in Fun(\V,\V')$  and $f\in Fun_\otimes(\F,\F')$ which commute with the structural maps $\imath,\imath'$ and $\imath^\otimes,\imath^{\prime \otimes}$ in the natural fashion. For simplicity, we assume that this means strict commutation. In general, these should be 2--commuting, see \cite{feynman}. Given such a morphisms the functor $f^*:\F'-\opcat_\C\to   \F'-\opcat_\C$ is simply given by precomposing $\O\mapsto f\circ \O$.

The push--forward is defined to be the  left Kan extension $Lan_f\O$. It has a similar formula as \eqref{kandefeq}. One could also write $f_!$ for this push--forward. Thinking geometrically $f_*$ is more appropriate.

We will reserve $f_!$ for the right Kan extension, which need not exist and need not preserve strong symmetric monoidality.  However, when it does it provides an extension by 0 and hence a triple of adjoint functors $(f_\ast, f^\ast,f_!)$.  This situation is characterized in \cite{Ward17} which also gives a generalization of $f_!$ and its left adjoint in those cases where the right Kan extension does not preserve strong symmetry.

\subsection{Examples}

\subsubsection{Tautological example}
$(\V,\V^{\otimes},\jmath)$. Due to the universal property of the free symmetric monoidal category, we have $\smodcat_\C\simeq\opcat_\C$.

{\sc Example:} If $\V=\underline G$, that is $\V$ only has one object, we recover the motivating example of group theory in the Warm Up. For a functor $\underline f: \underline G\to \underline H$
we have the functor $\underline{f}^{\otimes}$ and the pair $(\underline{f},\underline{f}^{\otimes})$ gives a morphism of Feynman categories. Pull--back becomes  restriction and  push--forward becomes
induction under the equivalence $\smodcat_\C\simeq\opcat_\C$.

Given any Feynman category $(\V,\F,\imath)$ there is always the morphism of Feynman categories given by $\imath$ and $\imath^\otimes$:
$(\V,\V^{\otimes},\jmath)\to (\V,\F,\imath)$ and the push--forward along it is the free functor $F$.

%

\subsubsection{Finite sets and surjections: $\F=Surj$, $\V=\trivial$} An instructive example for the hereditary  condition (ii) is the following.
As above let
 $Surj$  the category of finite sets and surjection with disjoint union $\amalg$ as monoidal structure  and let
$\trivial$  the trivial category with one object  $*$ and one morphism $id_*$.

 $\trivial^{\otimes}$ is equivalent to the category $\alephnsym$, where we think
 $\underline n=\{1,\dots,n\}=\{1\}\amalg\dots \amalg\{1\}$, $1=\imath (*)$.
  This identification ensures condition (i): indeed $\trivial^{\otimes}\simeq Iso(\Surj)$.

 Condition (ii) is more interesting. The objects of $(\F\downarrow \V$)  are the surjections $S\twoheadrightarrow\imath(*)$. Now consider an
arbitrary morphism of $\Surj$ that is a surjection $f:S \twoheadrightarrow T$ and pick an identification
 $T\simeq \{1,\dots,n\}$, where $n=|T|$. Then we can decompose the morphism $f$ as follows.
\begin{equation}
\xymatrix{
S\ar[rr]^{f}\ar[d]^{\simeq}&&T\ar[d]^{\simeq}\\
\amalg_{i=1}^{|T|} f^{-1}(i)\ar[rr]^{\amalg f|_{f^{-1}(i)}}&&\amalg_{i=1}^{|T|} \imath(*)
}
\end{equation}
Notice that both conditions (a) and (b) of \S\ref{detaildefpar} hold for these diagrams. This is because the fibers of the morphisms are well defined.
Condition (iii) is immediate. So indeed $\Surj=(Surj, \trivial, \imath)$ is a Feynman category.

 $\trivial$-$\smodcat_C$ is just $Obj(\C)$ and $Surj$-$\opcat$ are commutative and associative algebra objects or monoids in $\C$ as discussed in the Warm Up.
 The commutativity follows from the fact that if $\pi$ is the surjection $\underline 2\to \underline 1$, as above, and $\tau_{12}$ is the permutation of $1$ and $2$ in $\underline{2}=\{1,2\}$, which is also the commutativity constraint, then $\pi\circ \tau_{12}=\pi$ .

The functor $G$ forgets the algebra structure and the functor $F$ associates to every object $X$ in $\C$ the symmetric tensor algebra of $X$ in $\C$. In general, the commutativity constraints define what ``symmetric tensors'' means.

The monadicity can be read as in the Warm Up. Being an algebra over $GF$ means that there is one morphism for each symmetric tensor power $A^{\odot n}\to A$, that  on elements is given by $a_1\odot \dots \odot a_n\to a_1\dots  a_n$. This is equivalent to defining a commutative algebra structure.

The length of the morphisms is always non--negative and only isomorphisms have length $0$.

\subsubsection{Similar examples}
There are more  examples in which $\V$ is trivial and $V^{\otimes}\simeq \alephsym$.

Let $\F=Inj$ the category of finite sets and injections. This is  a Feynman category in which
 all the morphisms have non--positive length, with the isomorphisms being the only morphisms of length $0$.
If we regard $(\F\downarrow \V)$, we see that the injection $i:\emptyset\to \imath(\ast)$ is a  non--isomorphism, where $\emptyset=\Eins$ is the monoidal unit with respect to $\amalg$. By basic set theory, any other injection can be written as $id\amalg \dots \amalg id \amalg i \dots \amalg i$ followed by a permutation. This gives the decomposition for axiom (ii). The other two axioms are straightforward.

Using both injections and surjections, that is $\F=FinSet$, the category of finite sets and all set maps, we get the Feynman category $\Finset=(\trivial,FinSet,\imath)$.

\subsubsection{Skeletal versions: Biased vs.\ Unbiased. }
Notice that the skeletal versions of Feynman categories do give different $\opers$, although the categories $\opcat$ are equivalent. This is sometimes distinguished by calling the skeletal definition biased vs.\ the general set definition which is called unbiased. This terminology is prevalent in the graph based examples, see \S\ref{graphsubsec} and \S\ref{zoopar}.

\subsubsection{FI--modules and crossed simplicial groups, and free monoidal Feynman category}
We can regard the skeletal versions of the $\F$ above. For $sk(Inj)$ the ordinary functors $Fun(sk(Inj),\C)$ are exactly the FI--modules of \cite{Farb}.
Similarly, for $\Delta_+S$ the augmented crossed simplicial group, $Fun(\Delta_+S,\C)$ are augmented symmetric simplicial sets in $\C$.

In order to pass to symmetric monoidal functors, that is $\opcat$, one can use a free monoidal construction $\F^{\boxtimes}$. This associates to any Feynman category $\F$  a new Feynman category $\F^{\boxtimes}$ for which $\F^{\boxtimes}$-$\opcat_{\C}$ is equivalent
to the category of functors (not necessarily monoidal)
$Fun(\F,\C)$, see \S\ref{freemsec}.

\subsubsection{Ordered examples} As in the warm up, we can consider $\V=\trivial$, but look at ordered finite sets $\Finset_{ord}$ with morphisms being surjections/injections/all set morphisms. In this case the automorphisms of a set act transitively on all orders. For surjections we obtain not necessarily commutative algebras in $\C$ as $\opers$.

\subsection{Units}
Adding units corresponds to adding a morphisms $u:\emptyset\to\imath(\ast)$ and the modding out by the unit constraint $\pi\circ id_1 \otimes u=id_1$.
An $\oper$ $\O$ will take $u$ to $\eta=\O(u):\Eins\to A=\O(1)$.

\subsection{Graph Examples}
\label{graphsubsec}
\subsubsection{$\opcat$}
There are many examples based on graphs, which are explained in detail in the next section \S\ref{zoopar}. Here the graphs we are talking about are not objects of $\F$, but are
part of the underlying structure of the morphisms, which is why they are called ghost graphs. The maps themselves are morphisms between aggregates (collections) of corollas.  Recall that a corolla is a graph with one vertex and no edges, only tails. These morphisms come from an ambient category of graphs and morphisms of graphs. In this way, we obtain several Feynman categories by restricting the morphisms to those morphisms whose underlying graphs satisfy certain (hereditary) conditions.
The $\opcat$ will then yield types of operads or operad like objects.
As a preview:

\vskip\parskip

\begin{tabular}{l|l}
$\opcat$&Graph, i.e.\ underlying ghost graphs are of the form\\
\hline
Operads&  rooted trees\\
Cyclic operads&  trees\\
Modular operads& connected graphs (add genus marking)\\
PROPs&directed graphs (and input output marking)\\
NC modular operad& graphs (and genus marking)\\
Broadhurst-Connes&1-PI graphs\\
-Kreimer&\\
\dots&\dots\\
\end{tabular}

\vskip\parskip

Here the last entry is a new class. There are further decorations, which yield the Hopf algebras appearing in \cite{brown}, see \cite{deco}.

\subsubsection{Non--$\Sigma$ Feynman categories. The augmented simplicial category}
If we use $\V=\trivial$ as before, we can see that  $\F=\Delta_+$ yields a Feynman category. Now the non--symmetric $\V^{\otimes}=\alephnsym$
and the analog of $Surj$ and $Inj$ will then be order--preserving surjections and injections. These are Joyal dual to each other and play a special role in the Hopf algebra considerations.

Another non--$\Sigma$ example comes from planar trees where $\V$ are rooted planar  corollas and all morphisms preserve the orders given in the plane. The $\fopcatc$ are then non--sigma operads. Notice that a skeleton of $\V$ is given by corollas, whose in flags are labelled $\{1,\dots,n\}$ in their order and these have no automorphisms.

\subsubsection{Dual notions: Co--operads, etc.}  In order to consider dual structure, such as co--operads, one simply considers $\fopcat_{\C^{op}}$.
Of course one can equivalently turn around the variance in the source and obtain the triple:
$\FF^{op}=(\V^{op},\F^{op},\imath^{op})$. Now $\V^{op}$ is still a groupoid and $\imath^{\op,\otimes}$ still induces an equivalence,
but $\F^{op}$ will satisfy the dual of (ii). At this stage, we thus choose not to consider $\FF^{op}$, but it does play a role in other constructions.

\subsection{Physics connection} The name Feynman category was chosen with  physics in mind. $\V$ are the interaction vertices and
the morphisms of $\F$ are  Feynman graphs. Usually one decorates these graphs by fields.

In this setup, the  categories $(\F\downarrow \ast)$ are the channels in the $S$ matrix.
The external lines are given by the target of the morphism. The comma/slice category over a given target is then a categorical version of the $S$--matrix.

The functors $\O\in \fopcatc$ are then the correlation functions.
The constructions of the Hopf algebras agrees with these identifications and leads to further questions about identifications of various techniques in
quantum field theory to this setup and vice--versa. What corresponds to algebras and plus construction, functors? Possible answers
could be accessible via Rota--Baxter equations and primitive elements \cite{KKreimer}.

\subsection{Constructions for Feynman categories}
There are several constructions which will be briefly discussed below.

\begin{enumerate}
\item Decoration $\Fdeco$: this allows to define non--Sigma and dihedral versions. It also yields all graph decorations needed for the zoo; see \S\ref{decopar}.
\item $+$ construction and it quotient $\FF^{hyp}$: This is used for twisted modular operad and twisted versions of any of the previous structures; see \S\ref{enrichedpar}.
\item  The free constructions $\FF^{\boxtimes}$, for which $\FF^{\boxtimes}$-$\opcat_\C=Fun(\F,\C)$, see \S\ref{enrichedpar}.
Used for the simplicial category, crossed simplicial groups and FI--algebras.
 \item The non--connected construction $\FF^{nc}$, whose $\F^{nc}$-$\opcat$ are equivalent to lax monoidal functors of $\F$, see \S\ref{enrichedpar}.
\item The Feynman category of universal operations on $\FF$--$\opcat$ ; see \S\ref{univmasterpar}.
\item Cobar/bar, Feynman transforms in analogy to algebras and (modular) operads; see \S\ref{univmasterpar}.
\item W--construction, which gives a topological cofibrant replacement; see \S\ref{modelpar}.
\item Bi- and Hopf algebras from Feynman categories; see \S\ref{Hopfpar}.
\end{enumerate}

\section{Graph based examples. Operads and all of the Zoo}
\label{zoopar}
In this section, we consider graph based examples of Feynman categories. These include operads, cyclic operads, modular operads, PROPs, properads, their wheeled and colored versions, operads with multiplication, operads with $A_{\infty}$ multiplications etc., see Table \ref{zootable} They all come from a standard example of a Feynman category called $\GG$ via decorations and restrictions \cite{feynman,deco}. The category $\GG$ is a subcategory of the category of graphs of Borisov-Manin \cite{BM} and decoration is a technical term explained in section \S\ref{decosec}.

{\it Caveat:} Although $\GG$ is obtained from a category whose objects are graphs, the objects of the Feynman category are rather boring graphs; they have no edges or loops.
The usual graphs that one is used to in operad theory appear as underlying (or ghost) graphs of morphisms defined in \cite{feynman}. These two levels should not be confused and differentiate our treatment from that of \cite{BM}.

\subsection{The Borisov-Manin category of graphs.}
We start out with a brief recollection of the category of graphs given in \cite{BM}
\begin{enumerate}
\item A graph $\G$ is a tuple $(F_\G,V_\G,\del_\G,\imath_\G)$ of flags $F_\G$, vertices $V_\G$, an incidence relation $\del_\G:F\to V$ and an involution $\imath:F_\G^{\circlearrowleft}$, $\imath_\G^2=id$
which exhibits that either  two flags, aka.\ half-edges are glued to an edge  in the case of an orbit of order 2,  or  a flag  is an unpaired half--edge, aka.\ a tail if its orbit is of order one.

\item A graph morphism $\phi:\G\to \G'$ is a triple $(\phi_V,\phi^F,\imath_{\phi})$, where $\phi_V:V_\G\to V_{\G'}$ is a surjection on vertices, $\phi^F:F_{\G'}\to F_\G$ is an injection and
$\imath_{\phi}:F_\G\setminus \phi^{F}(F_{\G'})^{\circlearrowleft}$ is a self--pairing ($\imath_\phi^2=id$ and there are no orbits of order $1$).  This pairs together flags that ``disappeared'' from $\F_\G$ to ghost edges.
\item These morphisms have to satisfy  obvious compatibilities, see \cite{BM} or \cite{feynman}. One of these is preservation of incidence
$\phi_V\circ \del_\G\circ \phi^F(f')=\del_{\G'}(f')$ and ghost edges are indeed contracted $\phi_V(f)=\phi_V(\imath_\phi(f))$.

\end{enumerate}

We will call an edge $\{f,\imath(f)\neq f\}$ with two vertices $(\del(f)\neq \del(\imath(f))$ a simple edge  and an edge with one vertex $(\del(f)= \del(\imath(f))$ a simple loop.

As objects, the corollas are of special interest. We will write $\ast_S=(S,\{\ast\},\del: S\twoheadrightarrow \{\ast\},id)$ for the corolla with vertex $\ast$ and flags $S$.
This also explains our notation for elements of $\V$ in general.

An essential new definition \cite{feynman} is that of a ghost graph of a morphism.
\begin{df}
The ghost graph (or underlying graph) of a morphisms $\phi=(\phi_V,\phi^F,\imath_{\phi})$ is the graph $\gh(\phi)=(V_{\G},F_{\G},\hat\imath_{\phi})$, where $\hat \imath_\phi$ is the extension of $\imath_\phi$ to all of $F_\G$ by the identity on $F_\G\setminus \phi^F (F_{\G'})$.
\end{df}

\begin{ex}
Typical examples are isomorphisms ---which only change the names of the labels---, forming of new edges, contraction of edges and mergers.
The latter are morphisms which identify vertices. These identifications are kept track of by $\phi_V$. Composing the forming of a new edge and then subsequently contracting
it, makes the two flags that form the edge ``disappear'' in the resulting graph. This is what $\imath_\phi$ keeps track of. The ``disappeared'' flags form a ghost edge and this is the only way that flags may ``disappear''.
The ghost graph says that the morphism factors through a sequence of edge formations and subsequent contractions, namely those edges in the ghost graph, see Figure \ref{ghostgraphfig}.
\end{ex}

\begin{rmk} As can be seen from these examples:
The ghost graph does not determine the morphism. All the information about isomorphisms and almost all information about mergers is forgotten when passing from a morphism to the underlying graph.

What the ghost graph does, however, is keep track of are edge/loop contractions and this can be used to restrict morphisms. Further information is provided by the connectivity of the ghost graph, especially when mapping to a corolla. In this case, we see that mergers have non-connected ghost graphs.
Likewise, if we know that there are no mergers, then each component of the ghost graph corresponds to a vertex $v\in V_{\G'}$.
\end{rmk}

\subsubsection{Composition of ghost graphs corresponds to insertion of graphs into vertices}
The operation of inserting a graph $\gh_v$ into a vertex $v$ of a graph $\gh_1$,  is well defined for a given identification of the tails of $\gh_v$ with the flags $F_v$ incident to $v$. The result is the graph $\gh_v\circ_v \gh_1$
whose vertex set is $V=V_{\gh_1}\setminus \{v\}\amalg V_{\gh_v}$, the flags   $F=F_{\gh_1} \amalg F_{\gh_v}\setminus tails(\gh_v)$ with $\imath$  given by the disjoint union and $\del$ given by the disjoint union and the identification of $F_v$ with the tails of $\gh_1$.

Consider two composable morphisms and their composition:
$$
\xymatrix{
X
\ar@(dr,dl)[rr] _{\phi_0}
\ar[r]^{\phi_2}&Y\ar[r]^{\phi_1}&Z\\
\\
}
$$
\vskip -.5cm

Now let $\gh_i$ be the associated graphs of $\phi_i$, $i=0,1,2$.
Decomposing, $Y= \amalg_{v\in V_Y} \ast_v$,
and decomposing $\phi_2$ as $\amalg_{v\in V}\phi_v$
one can calculate  \cite{feynman} that $\gh_{0}$ is given by inserting each of the $\gh_v$ into the vertices $v$ of $\gh_{\phi_1}=V$, which we write as $\amalg_v \gh_v\circ \gh_1$.

\begin{equation}
\gh(\phi_0)=\gh(\phi_2)\circ \gh(\phi_1)
\end{equation}
where the identification for the composition is given by $\phi_2^F$.
%
%
%

\subsubsection{Symmetric monoidal structure}
The category of graphs has a symmetric monoidal structure given by disjoint union. The unit is the empty graph $(\emptyset,\emptyset,id_\emptyset,id_\emptyset)$
where $id_\emptyset:\emptyset\to \emptyset$ is the unique morphism from the empty set to itself.

\subsection{The Feynman category $\GG=(\Crl,\Agg,\imath)$}

Let
$\Crl$ be the subgroupoid of corollas with isomorphisms and $\Agg$.
$\Agg$ the full subcategory whose objects are aggregates of corollas. An aggregate of corollas  is a graph without any edges  $\imath_\G=id$.
Any aggregate of corollas is a (possibly empty( disjoint union of corollas and vice--versa.
Including corollas into the aggregates as one vertex aggregates gives an inclusion $\imath:\Crl\to\Agg$.

\begin{prop}
$\GG=(\Crl,\Agg,\imath)$ is a Feynman category.
\end{prop}

In this example the one--comma generators $(\F\downarrow \V)$ are morphisms from an aggregate to a simple corolla $*_v$

\begin{proof}
Looking at the definition of morphisms
it follows that $\Crl^\otimes\simeq Iso(\Agg)$. Condition (iii) is clear. For condition (ii) let $\phi:\G\to \G'$.  We will write any such morphism this as a disjoint union of one--comma generators.

For $v\in V_{\G'}$ define $\G_v$ to be the restriction of $\G$ to the vertices mapping to $v$. That is $\G_v=(V_{\G,v}=\phi_V^{-1}(v),F_{\G,v}= \del_\G^{-1}(V_{\G,v}),\del_{\G,v}=id)$. We let $\phi_v:\G_v\to v_{F_v}$ be the restriction of $\phi$, where $v_{F_v}$ is the corolla with vertex $v$ and its incident flags $F_v=\del_{\G'}^{-1}(v)$.
 $\G=(\phi_V|_{V_{\G,v}},\phi^F|_{F_v},\imath_\phi|_{F_{\G,v}\setminus (\phi^F)^{-1}(F_v)})$. It then follows that $\G=\amalg_{v\in V_\G'} \G_v,\G'=\amalg_{v\in V_{\G'}} v_{F_v}$ and $\phi=\amalg_{v\in V_{\G'}}\phi_v$. This yields the decomposition. It is easy to check conditions (a) and (b).
\end{proof}

\begin{figure}
    \centering
    \includegraphics[width=.7\textwidth]{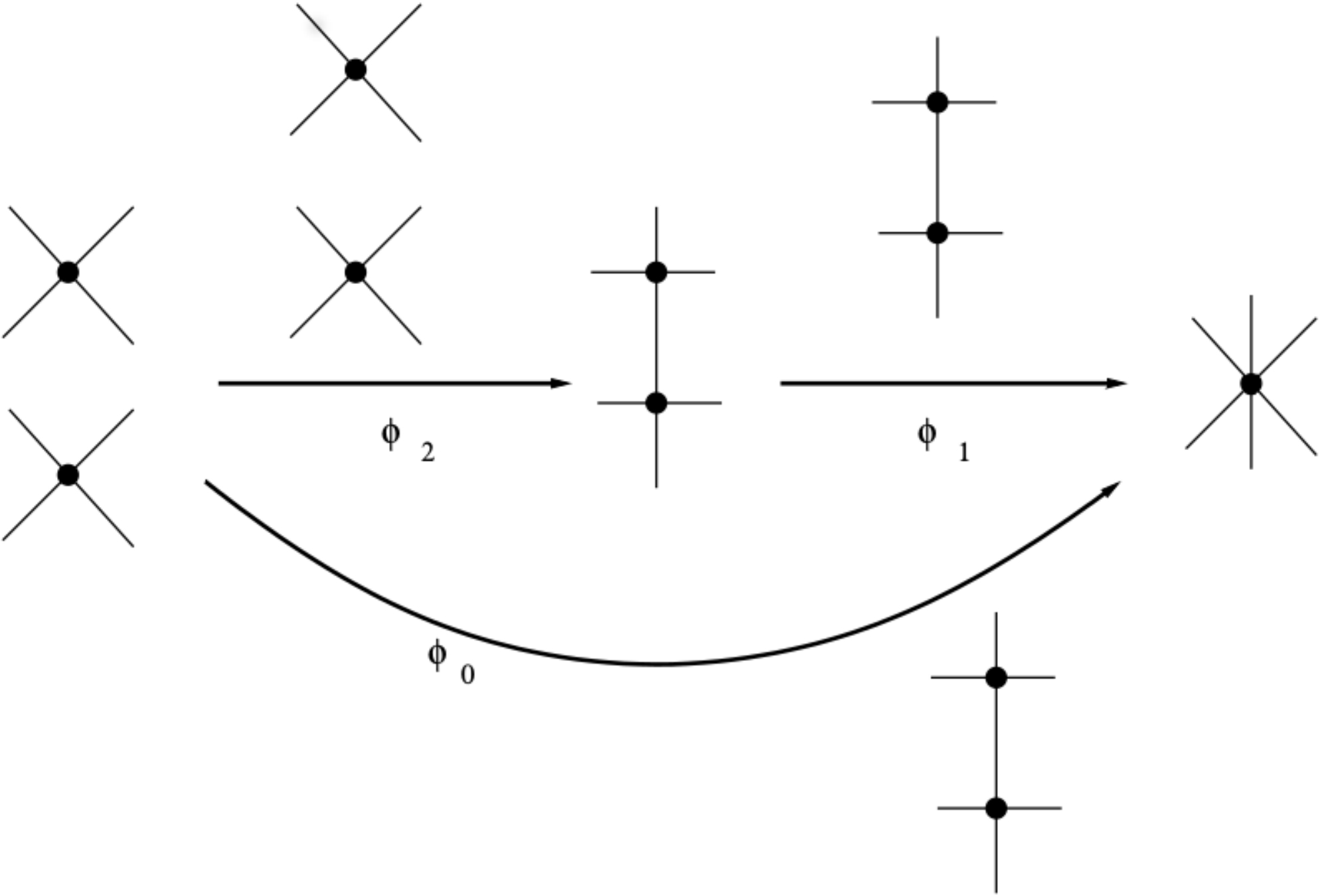}
    \caption{A composition of morphisms and the respective ghost graphs. The first morphism glues two flags to an edge, the second contracts an edge. The result is a morphism in $\Agg$.}
    \label{ghostgraphfig}
\end{figure}

Notice that forming an edge or a loop is not a morphism in $\Agg$. However the composition of the two morphisms, forming an edge or a loop {\em and the  subsequently contracting it} is a morphism in $\Agg$, see Figure \ref{ghostgraphfig}. One could call this a virtual or ghost edge contractions. For simplicity we will call these simply edge or loop contractions.

\subsubsection{Morphisms in $\Agg$}
\begin{figure}
    \centering
    \includegraphics[scale=.25]{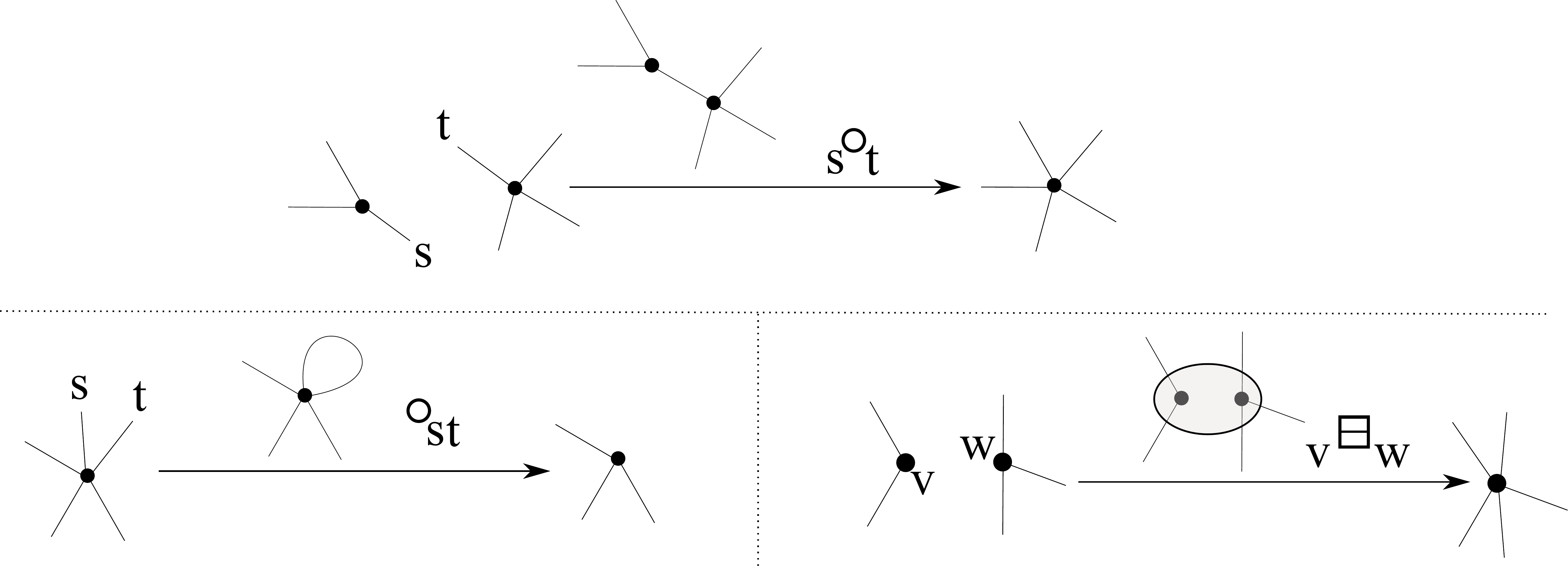}
    \caption{The three basic morphisms in $\GG$: an edge contraction (top), a loop contraction (left), and a merger (right). In the morphism, we give the ghost graph and label it by the standard notation. The shaded region is for illustration only, to indicate the merger.}
    \label{generatorfig}
\end{figure}

\begin{enumerate}
\item {\em Simple edge contraction.}  $\phi^F$ is the identity and the complement of the image $\phi^F$ is given by two flags $s,t$,
which form a unique ghost edge.
The two flags are not adjacent to the same
vertex and these two vertices are identified by $\phi_V$.  The ghost graph is obtained from the source aggregate by adding the edge $\{s,t\}$.
We will denote this by $\scirct$.

\item {\em Simple loop contraction.} As above, but  the two flags of the ghost edge
are adjacent to the same vertex. That is both $\phi_V$ and $\phi^F$ are identities. This is called a simple loop contraction.
 We will denote this by $\circ_{st}$.

\item {\em Simple merger.} This is a merger in which $\phi_V$ only identifies
two vertices $v$ and $w$. $\phi^F$ is an isomorphism.  Its degree is $0$ and the weight is $1$. The ghost graph is simply the source graph.
We will denote this by $\mge{v}{w}$.

\item {\em Isomorphism}. This is a relabelling preserving the incidence conditions. Here $\phi_V$ and $\phi_F$ are bijections. The ghost graph is the original graph.
\end{enumerate}
Typical examples of such morphisms are shown in Figure \ref{generatorfig}.

Actually any morphism is a composition of such morphisms \cite{feynman}. The relations  between these types of morphisms are spelled out below.
In order to make things canonical,
we will call a morphism pure  $\phi:\G\to \G'$, if $\phi^F=id$ when restricted to its image, and the vertices of $\G'$ are the fibers of $\phi_V$, that is $\phi_V(v)=\{w\in V_{\G}|\phi_V(w)=\phi_V(v)\}$.
With this terminology any morphism decomposes as
\begin{equation}
\label{bmdecompeq}
\phi=\sigma\circ\phi_m\circ \phi_c
\end{equation}
were $\phi_c$ is  a pure contraction, $\phi_m$ is a pure merger, and $\sigma$ is an isomorphism.

\subsubsection{Ghost graphs for $\Agg$}
In the case of morphism in $\Agg$, we can say more about the morphisms that have a fixed underlying ghost graph.
First, the source of a morphism  $\phi$ has the same vertices and flags as  its ghost graph  $\gh(\phi)$ and is hence completely determined.
If the ghost graph is connected, then  up to isomorphism the target is the vertex obtained from $\gh$  by contracting all edges.
If $\gh(\phi)$ is not connected, one needs the information of $\phi_V$  to obtain the target up to isomorphism. This is due to possible vertex mergers that are not recorded by the connected components of $\gh$. This information is encoded in a decomposition $\gh=\amalg_{v\in V} \gh_v$. The $\gh_v=\gh(\phi_v)$ are the ghost graphs of one--comma generators of the decomposition $\phi=\amalg_v \phi_v$.

Stated in another fashion:
in the decomposition \eqref{bmdecompeq}, $\gh(\phi)$ fixes $\phi_c$, the decomposition $\gh(\phi)=\amalg_v \gh_v$ fixes $\phi_m$.

\subsubsection{Relations}
\label{relpar}
All relations among morphisms in  $\GG$ are homogeneous in both weight and degree. We will not go into
the details here, since they follow directly from the description in the appendix. There are
the following types.

\begin{figure}
    \centering
 \includegraphics[scale=.25]{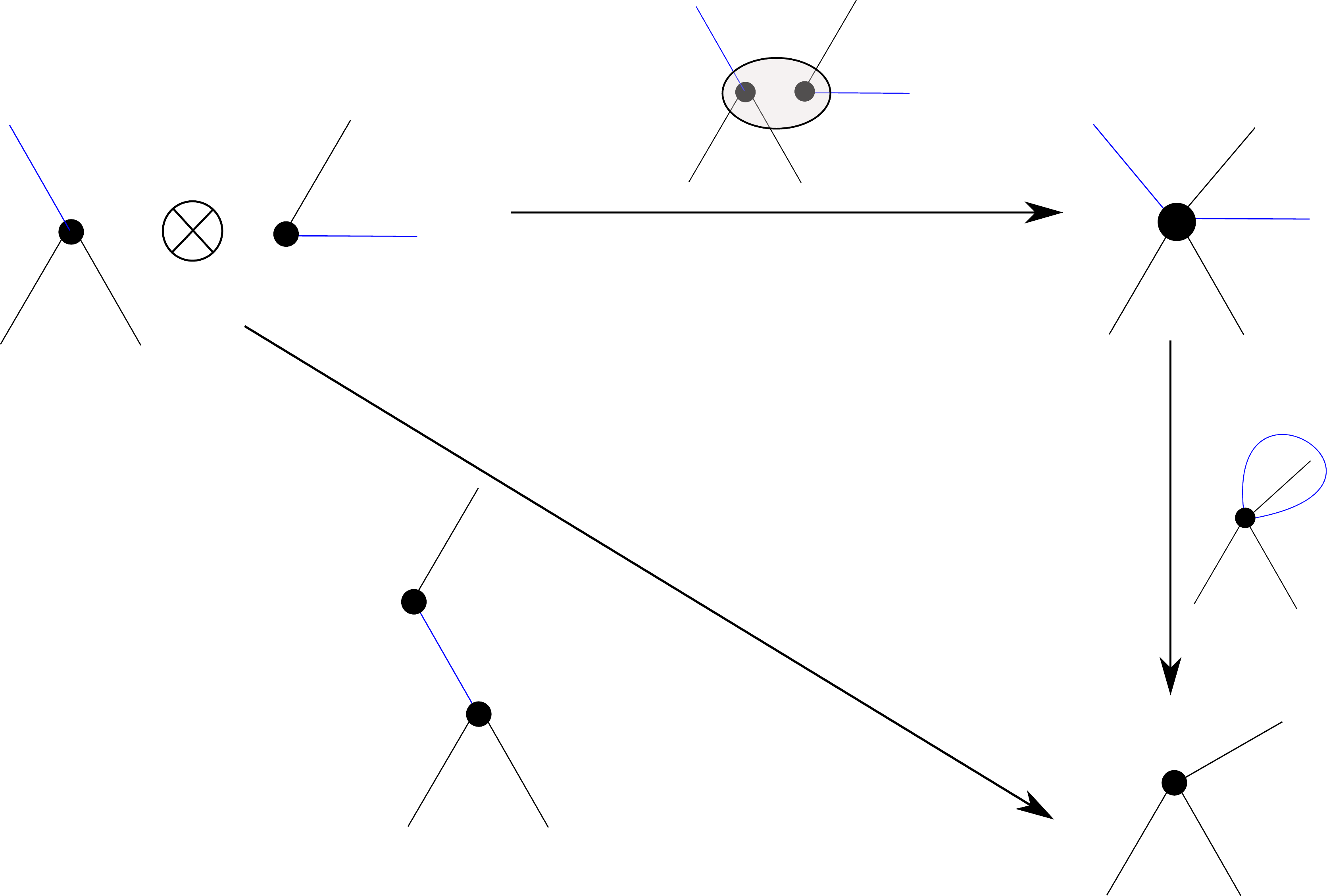}
    \caption{A triangle representing commutation between edge contraction and a merger followed by a loop contraction. The ghost graphs are shown. The shaded region is for illustrative purposes only, to indicate the merger.}
    \label{trianglefig}
\end{figure}

\begin{enumerate}
\item {\em Isomorphisms}. Isomorphisms commute with any $\phi$ in the following sense.
For any $\phi$ and any isomorphism $\sigma$ there are unique $\phi'$ and $\sigma'$ with $\gh(\phi\circ \sigma)=\gh(\phi')$
such
that
\begin{equation}
\phi\circ \sigma =\sigma'\circ \phi'
\end{equation}

\item {\em Simple edge/loop contractions}. All edge contractions commute in the following sense:
If two edges do not form a cycle, then the simple edge contractions commute on the nose
\begin{equation}\scirct \ccirc{s'}{t'} = \ccirc{s'}{t'} \scirct
\end{equation}
The same is true if one is a simple loop contraction and the other a simple edge contraction:
\begin{equation}\scirct\circ_{s't'}=\circ_{s't'}\scirct\end{equation}
If there are two edges forming a cycle, this means that
\begin{equation} \scirct \circ_{s't'}=\ccirc{s'}{t'} \circ_{st}\end{equation}
This is pictorially represented in Figure \ref{squarefig}.

\begin{figure}
    \centering
    \includegraphics[scale=.20]{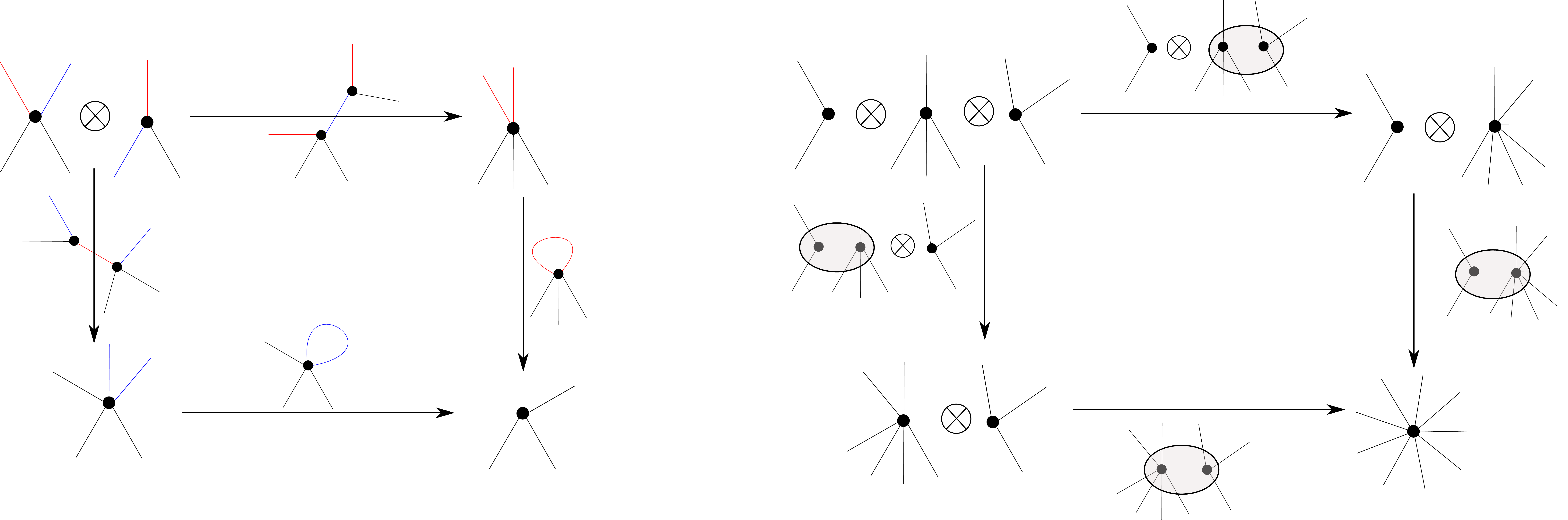}
    \caption{Squares representing commuting edge contractions and commuting mergers. The ghost graphs are shown. The shaded region is for illustrative purposes only, to indicate the merger.}
    \label{squarefig}
\end{figure}

\item {\em Simple mergers.} Mergers commute amongst themselves
\begin{equation}\mge{v}{w}\mge{v'}{w'}=\mge{v'}{w'}\mge{v}{w}\end{equation}
If $\{\del(s),\del(t)\}\neq\{v,w\}$ then
\begin{equation}
\label{mergereq}
\scirct\mge{v}{w} = \mge{v}{w}\scirct, \quad \circ_{st}\mge{v}{w}=\mge{v}{w}\circ_{st}
\end{equation}

If $\del(s)=v$ and $\del(t)=w$ then for a simple edge contraction, we have the following relation
\begin{equation}
\label{triangleeq}
\scirct = \circ_{st}\mge{v}{w}
\end{equation}
This is pictorially represented in Figure \ref{trianglefig}.
\end{enumerate}

\subsection{Examples based on $\GG$: morphisms have underlying graphs}
We are now  ready to present the zoo of operad--like structures in a structured way using the Feynman category $\GG$. The different Feynman categories will be obtained by decoration and restriction. Restriction often involves the underlying ghost graphs ---to be precise, the underlying ghost graphs of the one--comma generators. What one needs to check is that any such restriction is
stable under composition and the decorations compose, whence the term hereditary. For this it suffices to check compositions $X\to Y\to \imath(\ast)$. In other words, verify that $\amalg_v \gh_v\circ \gh$ satisfies a given restriction whenever $\gh$ and the $\gh_v$ are composable ghost graphs of one--comma generators satisfying this restriction. Likewise, one also has to define how the decorations compose and check that this gives an associative composition. The usual way is to induce the decoration on $\amalg_v \gh_v\circ \gh$ whenever the decorations on $\gh$ and the $\gh_v$ are given.
This can be done in the following cases (Table \ref{zootable}) in a straightforward fashion, see \cite{feynman} for details. For readers unfamiliar with some of these structures, the table may serve as a definition. We will discuss decorations, such as roots or directions in a more general fashion in \S\ref{decopar}.
For instance all these examples have {\em colored} versions by decorating the flags with colors.

We will say that $\FF$ is a Feynman category for a structure $X$ if $\F$-$\opcat_\C$ are the $X$--structures in $\C$.
E.g.\ $\operads$ is the Feymnan category for operads means that $\operads$-$\opcat_\C$ is the category of operads in $\C$.

\begin{table}[h]\begin{tabular}{llll}
$\FF$&Feynman category for&condition on ghost graphs $\gh_v$ and additional decoration&\\
\hline
$\operads$&(pseudo)--operads&rooted trees\\
$\operads_{May}$&May operads&rooted trees with levels\\
$\operads^{\neg\Sigma}$&non-Sigma operads &planar rooted trees\\
$\operads_{mult}$&operads with mult.&b/w rooted trees.\\
$\CCyclic$&cyclic operads&trees& \\
$\CCyclic^{\neg\Sigma}$&non--Sigma cyclic operads&planar trees& \\

$\GG$&unmarked nc modular operads& graphs \\
$\GG^{ctd}$&unmarked  modular operads&connected graphs \\
$\modular$&modular operads&connected + genus marking \\
$\modular^{nc,}$&nc modular operads &genus marking \\
$\dioperads$&dioperads&connected directed graphs w/o directed\\
&&loops or parallel edges\\
$\props$&PROPs&directed graphs w/o directed loops\\
$\properads$&properads&connected directed graphs \\
&&w/o directed loops\\
$\dioperads^{\circlearrowleft }$&wheeled dioperads&directed graphs w/o parallel edges \\
$\props^{\circlearrowleft,ctd}$& wheeled properads&connected directed graphs \\
$\props^{\circlearrowleft}$& wheeled props &directed graphs\\
$\FF_{1PI}$&1--PI algebras&1--PI connected graphs.
\end{tabular}

\caption{\label{zootable}List of Feynman categories with conditions and decorations on the graphs, yielding the zoo of examples}
\end{table}

New examples can also be constructed in this fashion. The first is the 1--PI (one particle irreducible) condition. A graph is 1PI if it is connected  furthermore even after remains connected after cutting any one edge the graph. There are more new examples of this type coming from quantum field theory and number theory, like the ones used in \cite{brown}, see \cite{GKT}.

\subsubsection{Push--forwards and pull--backs. Non--connected versions} There are obvious inclusion maps and forgetful maps between these categories. E.g.\ $\CCyclic\to \modular$, which assigns $g=0$ to each vertex. Here pull--back is the  restriction and push--forward is the modular envelope.
Looking at $\operads\to \props$, the root being ``out'', the push--forward is the PROP generated by and operad and the restriction is the operad contained in a PROP.
An examples that has been described by hand \cite{hoch2} is the PROP obtained from a modular operad. For this there is the morphism $\props\to \modular$,
which forgets the directions and adds genus $0$ to the vertices. Another is the inclusion $\modular\to \modular^{nc}$ which under push--forward gives the non--connected versions used for moduli spaces in \cite{KWZ,HVZ,schw}.

Analogously there is an inclusion $\FF\to \FF^{nc}$ for any of the candidates $\FF$ with connected graphs, where $\FF^{nc}$ allows non--connected graphs of the same type. Even more generally for and $\FF$ there is such a non--connected version $\FF^{nc}$ whose category $\opcat$ is equivalent to lax monoidal functors from $\FF$, see \S\ref{ncsec}.

\subsection{Details}

\subsubsection{Operad-lingo and notation: Composition along graphs, self gluing, non--self gluing and horizontal composition}

Let us unravel  the data involved in an $\O\in \F$-$\opcat$. Given a one--comma generator $\phi:X=\amalg_i \ast_{S_i}\to \ast_T$ we get a morphisms
$\O(\phi):\O(X)=\bigotimes_i \O(\ast_{S_i})\to \O(\ast_T)$. Here $X=s(\phi)$ is also the set of vertices of $\gh(\phi)$. If $\phi=\phi_c$ it is completely determined by its ghost graph and for pure contractions to corollas, which have connected ghost graphs,  we can set $\O(\gh(\phi)):=\O(\phi)$. This yields usual operad--like notations as follows. Define $\O(S):=\O(\ast_S)$. Then one can use the abbreviated notation $$\O(\gh):=\O(\gh(\phi)):\bigotimes \O(S_i)\to \O(T)$$ for the composition ``along
any connected graph $\gh$''.

For a simple edge contraction $\scirct:\ast_S\otimes \ast_T\to \ast_{(S\setminus s)\amalg (T\setminus t)}$ we get the standard non--self gluing pseudo operad compositions $\O(S)\otimes \O(T)\to  (S\setminus s)\amalg (T\setminus t)$, which is often denoted by $\scirct$ as well. In a similar maner, one obtains the May operations $\gamma$ for a rooted tree whose internal edges are all incident to the root.
A simple loop contraction $\circ_{s,s'}: \ast_{S}\to \ast_{S\setminus \{s,s\}}$ becomes the self gluing operation $\O(S)\to \O(S\setminus \{s,s'\})$; again by abuse of notation simply denoted $\circ_{s,s'}$.

If $\boxminus:\ast_S\amalg \ast_T\to \ast_{S\amalg T}$ is  a simple merger then in the usual
PROP notation this becomes the horizontal composition $\O(S)\otimes \O(T)\to \O(S\amalg T)$ usually also denoted by $\boxminus$.

Finally there are the isomorphisms. These are already incorporated into the $\V$--$\smodcat$ structure and not mentioned as structure operations in the operad--lingo. They are pushed into the underlying notion of $\SS$--module, or $\V$--$\smodcat$ in general, on which operads are built.
Thus by using \eqref{bmdecompeq} we can write any $\O(\phi)$ in the usual operad-lingo. The downside is that we {\em have to make this decomposition first}.

\subsubsection{Biased and unbiased versions}
Sending $S\to \ast_S$ provides an  equivalence from $\Finset$ to $\Crl$. We see that a skeleton of $\Crl$ is given by $\alephsym$.
Choosing $\V=\alephsym$, the $\V$--$\smodcat$ become $\alephsym$-modules.  Here usually one identifies $n$ with $\{0,1,\dots, n\}$ with $0$ indexing the root
if there is one present.

If we fix $Iso(\F)=\V^\otimes$ with $\V=\alephsym$, we obtain the biased notions of operads etc., that is objects $\O(n)$ with extra operations.
 Using $\V=\Finset$, we get $\O(S)$ with extra operations indexed by flags.

 If there is an extra decoration, then this is part of $\V$ and he set of vertices becomes bigger. An example is the genus marking in the modular operad case, so that we get $\O(n,g)$ or $\O(n,m)$ for
 Props, where $n$ are the incoming flags and $m$ are the outgoing flags in the biased version and $\O(S,g)$ and $\O(S,T)$ in the unbiased one.

 For instance, in the directed case a typical element of $\V$ is $\ast_{S,T}$ where
$S$ are the in--flags and $T$ are the out flags. Hence one obtains $\O(S,T)$ as for PROPs. Similarly if there is a genus marking a typical element is $\ast_{S,g}$ and hence in operad--lingo, we get $\O(S,g)$.

{\sc Variations}. If one is dealing with roots, often one uses the sets $n_+=\{0,\dots,n\}$ with the $0$ being the label of the root. An isomorphism must fix the roots, so that $Aut(\ast_{n_+})=\SS_n$. For operads, we then have the translation $\circ_i:=\ccirc{i}{0}$. In cyclic and modular operads, one commonly writes $\O((n))$  for $\O((n-1)_+)$ when using cyclic or modular operads, but does not insist that the maps are pointed, i.e.\ that the label $0$ is preserved, so that $Aut((n))=\SS_{n}$.

\subsubsection{A special case: PROP(erad)s vs.\ di-operads and wheeled versions}
\label{propsec}
PROPs and properads are a special case. Here the generators are not only the single edge contraction, but all multiparallel edge contractions.
In the graphs, parallel edges in the same direction are allowed. These cannot be factored into single edge constructions, so that there are generators
$\circ^k_{v,w}$ which simultaneously contract $k$ ghost edges of (necessarily) the same orientation between $v$ and $w$.

Allowing only the single edge contractions, one arrives at di--operads. Allowing wheels also allows to factor a multi--edge contraction and a single edge contraction followed by single loop contractions.
\subsubsection{Identities, multiplications etc as morphisms and decorations.}

We will briefly describe how to incorporate these operations.
Say, we want to add a ``unit'' as to get the Feynman category for unital operads. Recall that for and operad  $\O$ a unit is an element $\eta:\Eins_\C \to \O(1)$
which satisfies $u\circ_1 a=a=a\circ_i u$.

\begin{figure}[htb]
\begin{center} \includegraphics[scale=.4]{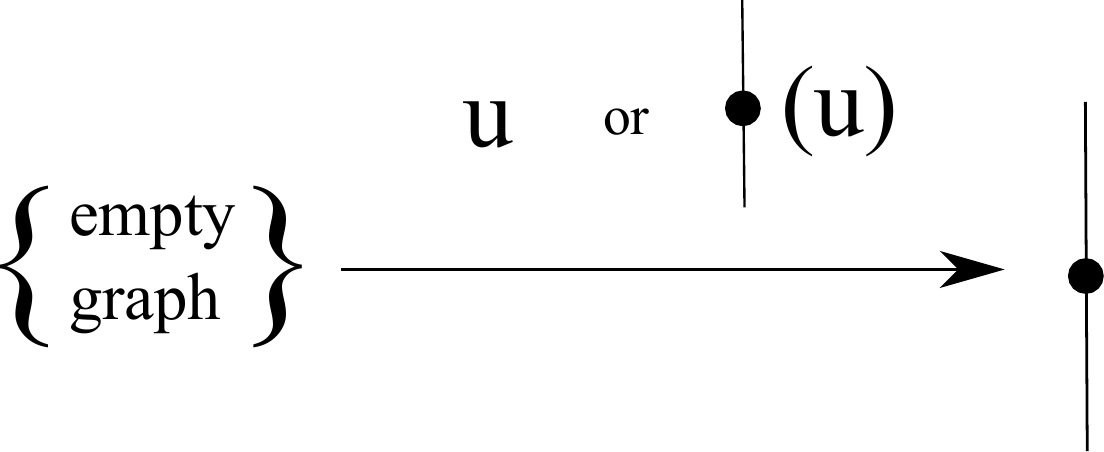} \end{center}
\caption{\label{unitfig}Graphically adding a morphism as marked binary vertex of the ghost graph}
\end{figure}

Since $\Eins_\C=\O(\Eins_\F)$, we adjoint a morphisms $u:\emptyset\to \ast_{1_+}$ to  the Feynman category for operads $\operads$ with source the empty graph. This can be graphically noted by putting a
$u$ on a binary vertex of a ghost tree, whenever we want to use the morphism $u$, as illustrated in Figure \ref{unitfig}.
This does not yet constitute putting in a unit, but rather asking for the data of an element in $\O(1)$.
This is actually what is needed in the case of the Hopf algebra of Connes and Kreimer \cite{GKT}, see also \S\ref{Hopfpar}.
In order to get a unit, we have to quotient by the relation given above. The simplest graphical way to do this is to remove all the vertices $u$ from the graph. Technically this is given by an equivalence relation. If one does this, one can create  a new ``degenerate graph'' consisting of a lone flag, which represents any tree whose vertices are all marked by $u$. This explains the notation of e.g.\ \cite{Markl}.

In this fashion, one sees that one gets an isomorphism of Feymnan categories between the Feynman category for unital May operads and that for unital operads, see \cite{feynman} for details.

Similarly, for multiplications one needs an extra morphism $\mu:\unit_\C\to \O(2)$. Consequently, one adjoins a morphism $\emptyset \to \ast_{2_+}$. In the graphical version, the (ghost) graphs will now have a possible decoration on 3--valent vertices by $\mu$.
This just gives a multiplication, one can then quotient out by the associativity equation. This amounts to graphs with black and white vertices, where black indicates an iteration of $\mu$. Here associativity induces an equivalence relation, which allows to contract all edges of any subtree of vertices marked solely by $\mu$.
A similar procedure adds the $\mu_n$ for $A_\infty$ multiplications as black vertices of arity $n$, see e.g.\ \cite{woods,KSchw,feynman}.

Furthermore all these kinds of extra morphisms can be collected and turned into a decoration in the technical sense. This is detailed in \cite{feynman}.

\subsection{Omnibus Theorems}
 For any of these, we have a general triple of graphs $\T=GF$. We immediately obtain a general theorem for all of the zoo and all new species of this kind;
 see also, \S\ref{decopar}. These give the usual three  ways of describing these objects  (a) via composition along graphs, (b) as algebras over a triple or (c) via generators and relations for the morphisms.

 \begin{thm}
 The biased and unbiased $\opcat_\C$ are equivalent. Moreover the $\F$--$\opcat_\C$ are equivalent to algebras over the relevant triple of graphs.
 \end{thm}

Notice the usual triples of graph, see e.g.\ \cite{MSS}, match up exactly with the triples above, when one considers the ghost graphs and their composition.
Moreover, the whole semi--simplicial structure of iterating the endofunctors, cf.\ \cite{MSS,GKmodular}, coincides as demonstrated in \cite{feynman}.

\begin{thm} Generators and relations description. All the examples have a generator and relations description. The generators always contain the isomorphisms, the edge contractions $\scirct$. If non--connected graphs are allowed, the morphisms include the mergers $\boxminus_{v,w}$ and if loops are allowed, then they contain the loop contractions. In the presence of decorations, these are restricted to respect the decorations (cf.\ \S\ref{decopar}). The relations are the ones given above.

If one adds additional morphisms with relations, these are  be included in the list.
\end{thm}
This can be formalized using Feynman categories indexed over another Feynman category, see \cite{feynman}.

{\sc Example:} For instance, when adding units, the morphism $u$ is a generator and the relations with $u$ are the unit relations.
This way, one can, for example, get the Feynman category for unital cyclic operads in all three definitions.

\begin{rmk}
In the PROP(erad) case, which is special,  the generators are not only the simple edge contraction, but multi--edge contractions, see \S\ref{propsec}.
\end{rmk}

\section{Decorating Feynman categories $\FF_{dec\O}$ }
\label{decopar}

Decorations can be made into a technical definition. The details for this section are in \cite{deco}.
The basic idea is that one can decorate a Feynman category by using elements of  $\F$--$\opcat$.
The reason this works is that in order to define a composition, one has to give a composition for the decorations, but this is precisely the data of an $\O\in\F$--$\opcat$. These decorations actually decorate the elements of $\V$. In the graph example above, this means that one can decorate vertices and  flags.

\subsection{Main Theorems}

The main constructive theorem is the following.

\begin{thm} Given an $\O\in \F$--$\opcat$, then there is a Feynman category $\F_{dec\O}$ which is indexed over $\F$. It objects are pairs $(X,dec\in \O(X))$ and $Hom_{\F_{dec\O}}((X,dec),(X',dec'))$ is the set of $\phi:X\to X'$, s.t. $\O(\phi):dec\to dec'$. \end{thm}

\begin{rmk}
This theorem also works in the enriched setting, where one considers enrichment over $\C$, confer \S \ref{enrichedpar}.
This construction works directly for Cartesian $\C$, and with modifications it also  works for the non--Cartesian case.
\end{rmk}

\begin{ex} All planar structures:
Non--sigma operads, cyclic non--Sigma operads, non--Sigma modular operads.
Here $\O$ is $\mathcal{A}ssoc$, $\mathcal{C}yc\mathcal{A}ssoc$, $\mathcal{M}od\mathcal{C}yc\mathcal{A}ssoc$.
These are actually all obtained by functoriality, see below.  This recovers e.g.\  that the modular envelope of $\mathcal{C}yc\mathcal{A}ssoc$ factors through non--Sigma modular operads  \cite{marklnonsigma}.

\end{ex}
%
\begin{thm}[Functoriality in $\FF$ and $\O$]
\label{decofunctthm}
Given a morphism of Feynman categories $f:\FF\to\FF'$ and a morphisms $\sigma:\O\to \P$.
There are commutative squares which are natural in $\O$
\begin{equation}
\xymatrix{\Fepair \ar[r]^{f^{\O}} \ar[d]_{forget} & \Fe'_{dec\, f_{\ast}(\O)} \ar[d]^{forget'} \\
\Fe \ar[r]^f & \Fe'}
\quad
\xymatrix{\Fepair \ar[r]^{\sigma_{dec}} \ar[d]_{f^{\O}} & \Fe_{dec\Po} \ar[d]^{f^{\Po}} \\
\Fe'_{decf_{\ast}(\O)} \ar[r]^{\sigma'_{dec}}& \Fe'_{decf_{\ast}(\Po)} }
\end{equation}
On the categories of monoidal functors to $\C$, we get the induced diagram of adjoint functors.

\vskip -5mm

\begin{equation}
\xymatrix{\Fpair\text{-}\opcat \ar@/^/[r]^{f^{\O}_*}\ar@/_/[d]_{forget_*} & \F'_{dec\,f_{\ast}(\O)}\text{-}\opcat   \ar@/^/[l]^{f^{\O}*}\ar@/^/[d]^{forget'_*} \\
\F\text{-}\opcat \ar@/^/[r]^{f_*}\ar@/_/[u]_{forget^*} & \F'\text{-}\opcat\ar@/^/[u]^{forget'^*}\ar@/^/[l]^{f^*} }
\end{equation}
\end{thm}

\subsection{Terminal objects and minimal extensions}

\begin{thm}
\label{decothm}
If $\final$ is a terminal object for $\fops$ and $forget:\Fdeco\to \F$ is the forgetful functor,
then $forget^*(\final)$ is a terminal object for $\Fdeco\text{-}\opcat$.
We have that $forget_*forget^*({\final})=\O$.
\end{thm}

\begin{df}
We call a morphism of Feynman categories $i:\FF\to \FF'$ a minimal extension over $\C$ if  $\FF$-$\opcat_{\C}$ has a
 a terminal/trivial functor $\final$ and $i_*{\final}$ is a terminal/trivial functor in $\FF'$-$\opcat_{\C}$.
 \end{df}
\begin{ex}
There are two examples that appear naturally. The first is $CycCom$ and
$ModCycCom$ for $\CCyclic \to \modular$ and the second is the decorated version
$\forget^*(CycAssoc)$ and $i^{\O}_*(\forget^*(CycAssoc))$.
\end{ex}
\begin{prop}
\label{minextprop}
If $f:\Fe \to \Fe'$ is a minimal extension over $\C$, then $f^{\O}:\Fepair \to \Fepairtwo$ is as well.
\end{prop}


\subsection{Example}

\subsubsection{Markl's Non-$\Sigma$ modular (see also \cite{KP})}
 \begin{equation}
\label{modulardiag}
\xymatrix{\F_{dec\, CycAssoc}=\CCyclic^{\neg \Sigma} \ar[r]^{i^{CycAssoc}} \ar[d]_{forget} & \modular_{dec\, i_{\ast}(CycAssoc)} =\modular^{\neg \Sigma}\ar[d]^{forget} \\
\CCyclic \ar[r]^i & \modular}
\end{equation}

\begin{enumerate}
\item The commutative square exists simply by Theorem \ref{decofunctthm}.
\item  On the left side,  if $*_C$ is final for $\CCyclic$ and hence $forget^*(*_C)=\underline{*}_C$
is final for $\CCyclic^{\neg \Sigma}$ . The pushforward $forget_*(\underline{*}_C)=CycAssoc$.
\item On the right side,  if $*_M$ is final for $\modular$ and hence $forget^*(*_M)=\underline{*}_M$
is final for $\modular^{\neg \Sigma}$. The pushforward $forget_*(\underline{*}_M)=ModAssoc$.
\item The inclusion $i$ is a minimal extension.
This is a fact explained by basic topology.
Namely gluing together polygons in their orientation by gluing edges pairwise yields all closed oriented surfaces, see e.g.\ \cite{munkres}.
\item  Hence $i^{CycAssoc}$ is also a minimal extension.
which explains why indeed the pushforward of the terminal $\oper$ is up to that point  still terminal.
It also reflects the fact that not gluing all edges pairwise, but preserving orientation, does yield all surfaces with boundary.
\end{enumerate}

\subsection{Examples on $\GG$ with extra decorations, non--sigma, colored versions, etc}

\label{decosec}

We now give the details on how to understand the decorations in \S\ref{zoopar} as decorations in the technical sense.
Decoration and restriction allows to generate the whole zoo and even new species. Examples of the needed decorations are listed in Table
\ref{decotable}.
\begin{table}[h]
\begin{tabular}{llll}
$\FFdeco$&Feynman category for&decorating $\O$&restriction\\
\hline
$\FF^{dir}$&directed version&$\Z/2\Z$ set&edges contain one input\\
&&& and one output flag\\
$\FF^{rooted}$&root&$\Z/2\Z$ set&vertices have one output flag.\\
$\FF^{genus}$&genus marked&$\N$&\\
$\FF^{c-col}$&colored version&$c$ set&edges contain flags\\
&&&  of same color\\
$\operads^{\neg\Sigma}$&non-Sigma-operads&$Assoc$&\\
$\CCyclic^{\neg\Sigma}$&non-Sigma-cyclic operads&$CycAssoc$&\\
$\modular^{\neg\Sigma}$&non--Signa-modular&$ModAssoc$&\\
$\CCyclic^{dihed}$&dihedral&$Dihed$&\\
$\modular^{dihed}$&dihedral modular&$ModDihed$&
\end{tabular}
\caption{\label{decotable}List of decorated Feynman categories with decorating $\O$ and possible restriction. $\FF$ stands for an example based on $\GG$ in the list or more generally indexed over $\GG$ (see \cite{feynman}). }
\end{table}

\subsubsection{Flag labelling, colors, direction and roots as a decoration}
Recall that $*_S$ is the one vertex graph with flags labelled by $S$ and these are the objects of $\V=\Crl$ for $\GG$.
For any set $X$ introduce the following $\GG$-$\oper$: $X(*_S)=X^{S}$. The compositions are simply given by restricting to the target flags.

Now let the set $X$ have an involution $\bar{}:X\to X$.
Then a natural subcategory $\FF_{decX}^{dir}$ of $\GG_{decX}$ is given by the wide subcategory, whose morphisms additionally satisfy that only flags marked by
elements $x$ and $\bar x$ are glued and then contracted; viz\ $\imath_\phi$ only pairs flags of marked $x$ with edges marked by $\bar x$.  That is the underlying ghost graph has edges whose two flags are labelled accordingly.
 In the notation of  graphs:  $X(f)=\overline{\imath_{\phi}(f)}$.

 If $X$ is pointed by $x_0$, there is the subcategory of $\GG_{decX}$ whose objects are those generated by $*_S$ with exactly one flag labelled by $x_0$
 and where the restriction on graphs is that for the underlying graph additionally, each edge has one flag labelled by $x_0$.

Now if $X=\Z/2\Z=\{0,1\}$ with the involution $\bar 0=1$,  we can call $0$ ``out'' and $1$ ``in''.  As a result, we obtain the category of directed graphs $\GG_{dec\Z/2Z}$.
Furthermore, if $0$ is the distinguished element, we get the rooted version. This explains the relevant
examples  Table \ref{decotable}.

More generally, in quantum field theory the involution sends a field to its anti--field and this is what decorates the lines or propagators in a Feynman graph.

\subsubsection{Genus decoration}
Let $\N$ be the $\GG$-$\oper$ which on objects of $\V$ has constant value the natural numbers $\N(*_S)=\mathbf{N_0}$. On morphisms $\N$ is defined to behave like the genus marking.
That is for $\phi:X\to *_S$, we define $\N(\phi):\N(X)=\mathbf{N_0}^{|X|}\to \mathbf{N_0}=\N(*_S)$
 as the concatenation $\mathbf{N_0}^{|X|}\stackrel{\sum}{\to} \mathbf{N_0}\stackrel{+\bar\gamma(\phi)}{\to}\mathbf {N}_0$  where $\bar\gamma(\phi)$ equals one minus the Euler characteristic of the graph underlying $\phi$. If this graph is connected this is just
first Betti number also sometimes called the genus. This coincides with the description in \cite{feynman}, Appendix A.
Hence, if $\FF$ is a subcategory of $\GG$, then the genus marked version  is just $\FF_{dec\N}$.
Examples are listed in Table \ref{decotable}.

\subsubsection{Assoc-decorated, aka.\ Non--Sigma, aka.\ non--planar}
Likewise, we can regard the cyclic associative operad, $CycAssoc$.
The pull back of $CycAssoc$ under $forget:\operads\to\CCyclic$ is the associative operad $Assoc$.
Now $\operads_{dec\, Assoc}=\operads^{\neg\Sigma}$ is the Feynman category for non--Sigma operads. Indeed, the
elements of $Assoc(*_s)$ are the linear orders on $S$, which means that we are dealing with planar corollas as objects.
Likewise, for the morphisms the condition that $\phi(a_X)=a_Y$ means that the trees are also planar.
The story for cyclic operads is similar $\CCyclic_{dec CycAssoc}=\CCyclic^{\neg\Sigma}$.

Things are more interesting in the modular case. In this case, we have
$ModAssoc:=i_*(CycAssoc)$ as a possible decoration and we get the decorated Feynman category $\modular^{\neg\Sigma}:=\modular_{dec \, ModAssoc}$..
Indeed using this decoration,  we recover the definition of \cite{marklnonsigma} of non--sigma modular operads, which is the special case of a brane--labelled  c/o system, with trivial closed part and only one brane color   \cite{KP}[Appendix A.6]; see also \cite{KLP},
the appendix of \cite{postnikov} and \cite{marklnonsigma} for details about the correspondence between stable or almost ribbon graphs and surfaces.

Here we can understand these constructions in a more general framework.
First, the diagram considered in   \cite{marklnonsigma} is exactly a diagram of Theorem \ref{decofunctthm}. Then the fact that  the non--Sigma modular envelope
of $CycAssoc$ is terminal is obvious from Theorem \ref{decothm} and Proposition \ref{minextprop}. The key observations
are that the terminal object of $\CCyclic^{\neg\Sigma}$  pushed forward is indeed $CycAssoc$ and that $ModAssoc$  is the pushforward
of the terminal object of  $\modular^{\neg\Sigma}$.
 Notice $CycAssoc$ is not a modular operad, so it is not a valid decoration for $\modular$. This is reflected in the treatments of \cite{KP,marklnonsigma}. We see that we do get a planar aka.\ non--Sigma version by pushing forward $Assoc$.

\subsection{Kontsevich's three geometries}
In this framework, one can also understand Kontsevich's three geometries \cite{kontsevichthree} as follows.

\subsubsection{Com, or trivially decorated}
The operad $CycCom$, the operad for cyclic commutative algebras, is the terminal/trivial object in $\CCyclic$-$\opcat$.
Thus by Theorem \ref{decothm},
we have that $\operads_{decCom}=\operads$. The analogous statement holds for $\CCyclic$. Indeed,
there is a forgetful functor $\operads\to\CCyclic$ and the pull--back of  $CycCom$ is $Com$ and hence
$\CCyclic_{decCycCom}=\CCyclic$. Finally using the inclusion $i:\CCyclic\to \modular$ means that the modular envelope
$i_*(Com)$ is a modular operad. Tracing around the trivially decorated diagram, we see that
this is again a terminal/trivial operad. Indeed this is the content of Proposition \ref{minextprop}.

\subsubsection{Lie, etc. or graph complexes} For this we actually need the enriched version.

One of the most interesting generalizations is that of Lie or
in general of Kontsevich graph complexes. Here notice that $Assoc,Com$ and $Lie$ are all three cyclic operads, so that they all can be used to decorate the Feynman category for cyclic operads. For $Lie$ it is important that we can also work over $k$--Vect.
Thus, answering a question of Willwacher \cite{Willwacherprivate}, indeed there is a Feynman category for the Lie case.

To go to the case of graph complexes, one needs to first shift to the odd situation and then  take colimits as described in detail in \cite{feynman}, see especially section 6.9 of {\it loc.\ cit.}.

\subsection{Further applications}
Further  forthcoming applications will be
\begin{enumerate}
\renewcommand{\theenumi}{\roman{enumi}}
\item Infinity versions of the Assoc, Com and Lie and their transformations.
\item New decorated interpretation of moduli space operations generalizing those of \cite{hoch1,hoch2}.
\item The new Stolz--Teichner--Dwyer setup for twisted field theories.
\item Kontsevich's graph complexes.
\item Actions of the Grothendieck--Teichm\"uller group.

\end{enumerate}

%
%
%
%
%
%
%
%

\section{Enrichment, algebras, odd versions and further constructions}
\label{enrichedpar}
\subsection{Enriched versions, plus construction, and algebras over $\FF$--$\opcat$.  Overview and examples}
There are several reasons why one would like to consider enriched versions of Feynman categories. They are necessary to define the transforms and resolutions.
Here it is necessary to introduce signs or anti--commuting morphisms. They are also natural from an algebra over operads point of view.
We will start with this construction.

\subsubsection{The Feynman category for an algebra over an operad}
\label{algebrasec}
Recall that an algebra over an operad $\O$ in $\C$ is an object $A$ and a morphism of operads $\rho:\O\to End(A)$. For this to make sense,
one assumes that $\C$ is closed monoidal. Then $\CalE nd(A)(n)=\underline {Hom}(A^{\otimes n},A)$. One can simply think of $\C=\Vect$ or $\Set$. Substitutions then give the operad structure.

{\sc Algebras as natural transformations.} Generally, given a reference target $\FF$-$\oper$ $\CalE$, then for another $\O\in \fopcatc$ we  define  an $\O$--algebra
relative to $\CalE$ as a natural transformation of functors $\rho:\O\to \CalE$.

Indeed,  for instance in the operad case with $\CalE=\CalE nd$, we obtain $\rho(n):\O(n)\to Hom(A^{\otimes n},A)$ which commute with compositions.

{\sc Algebras over operads as functors.}
We will start with the operad case. Given a May operad $\O$, we will construct a Feymnan category $\FF_\O$ whose $\opers$ are algebras over $\O$.
The data we have to encode are $A\in C$ and $\rho(n):\O(n)\to Hom(A^{\otimes n},A)$.
Now if we take $\V_\O=\trivial$ and $Iso(\F_\O)=\alephsym$, then we see that a strict symmetric monoidal functor $\rho:\alephsym\to \C$  will send $n$ to $A^{\otimes n}$ and
the $\sigma\in Aut(n)=\SS_n$ to the  permutations of the factors of $A^{\otimes n}$.

We now add more morphisms. A morphisms from $\phi:n\to 1$ will be sent to a morphism $\rho(\phi):Hom(A^{\otimes n},A)$. Thus, we set the one--comma generators as $\O(n)=:Hom_{\F_\O}(n,1)$. This fixes data of the $\rho(n)$ is  and vice--versa. Notice that when adding in these morphisms, $\O(n)$ is  ---and has to be--- an $\SS_n$--module to fix the pre--composition with the isomorphisms $Aut(n)$.

 Here we assume that we can also work with enriched categories. In particular, we need to be enriched over $\C$ if $\O$ is an operad in $\C$, see details below.

With these one--comma generators, due to condition (ii), we get that
$Hom_{\F_\O}(n,m)=\bigotimes_{(n_1,\dots,n_m):\sum n_i=n} \O(n_1)\odo \O(n_m)$. Here $\bigoplus$ is the colimit, which we assume to exist.
There is more data. In order to compose
$Hom_{\F_\O}(m,1)\otimes Hom_{\F_\O}(n,m)\to  Hom_{\F_\O}(n,1)$, we need morphisms
\begin{equation}
\gamma_{n_1,\dots,n_k}: \O(m) \otimes \O(n_1)\odo \O(n_m)\to \O(n) \quad n=\sum n_i
\end{equation}
These have to be compatible with the isomorphisms.
This data is the composition of a May operad and vice--versa defines a category structure on $\F_\O$.

This category has a special structure, namely that
\begin{equation}
Hom_{\F_\O}(n,m)=\bigoplus_{\phi:n\twoheadrightarrow m} \O(\phi) \text{ where } \O(\phi)=\bigotimes_{i\in m}\O(f^{-1}(i))
\end{equation}


{\sc Caveats:} In order to obtain a Feynman category, we will need to define what an enriched Feynman category over $\C$ is. This is straightforward if $\C$ is Cartesian. In the non--Cartesian case, we have to be a bit more careful, see below.
There we will see that the isomorphism condition will dictate that $\O(1)$ has only $\id$, that is a copy of $\unit_\C$ corresponding to $id$
as the ``invertible element''.
Also, the relevant notion is that of a Feynman category indexed enriched over another Feynman category.
In our example, we are indexed enriched over a skeleton of $\Surj$.

Clearing these up leads to the theorem:

\begin{thm}
The category of  Feynman categories  enriched over $\CalE$ indexed over $\Surj$ is equivalent to the category of operads (with the only iso in $\O(1)$ being the identity) in $\CalE$
with the correspondence given by $O(n)=Hom(\underline{n},\underline{1})$. The $\opcat$ are now algebras over the underlying operad.
\end{thm}

\begin{rmk}
We can also deal with algebras over operads which have isomorphismS in $\O(1)$ by enlarging $\V$. For this one needs a splitting $\O(1)=\O(1)^{iso}\oplus \bar\O(1)$, where no element of $\bar \O(1)$ is invertible and $\O(1)^{iso}=\bigoplus_{g\in G}\unit_\C$ for an index group $G$ is the free algebra on $G$. Then we enlarge $\V$ by letting $1$ have isomorphisms $G$. The construction is then analogous to the one above and that of $K$--algebras \cite{GKmodular}.
Another way is to use lax moniodal functors, see \cite{feynman}.
\end{rmk}

\subsubsection{General situation for algebras: Plus construction}
There is a "+" construction, not unlike that for polynomial monads \cite{monads}, that produces a new Feynman category out of an old one.
Inverting morphisms stemming from isomorphisms one obtains $\FF^{hyp}$ and there  is a further reduction to an equivalent category $\FF^{hyp,rd}$.
Details will be provided below.

The main theorem is that enrichments of $\FF$ are  in 1--1 correspondence with $\FF^{hyp}$--$\opcat$.

\begin{ex}
 $\modular^{hyp}=\FF_{hyper}$, the Feynman category for hyper--operads as defined by \cite{GKmodular}, whence the name.
 $\Surj^+=\FF_{Mayoperads}$,  $\FF_{surj}^{hyp,rd}=\operads_0$, the category for operads whose $\O(1)$, has only (multiples of) $id$ as an invertible element.
 $\FF_{triv}^+=\Surj$, $\FF_{triv}^{hyp,rd}=\FF_{triv}$.
\end{ex}

\begin{df}
Let $\FF$ be a Feynman category and $\FF^{hyp,rd}$ its reduced hyper category, $\O$ an
$\FF^{hyp,rd}$-$\oper$ and $\D_{\O}$ the corresponding enrichment functor. Then we define an $\O$-algebra to be a $\FF_{\D_{\O}}$-$\oper$.
\end{df}

\subsubsection{Odd Feynman categories over graphs}
\label{oddgraphpar}
In the case of underlying graphs for morphisms, odd usually means that edges get degree $1$, that is we use a Kozsul sign with that degree.
In particular, in these discussions, one is augmented over $\Ab$, the category of Abelian groups.
Then there is an indexed enriched version of the Feynman categories.
In order to write this down, one needs an ordered presentation.

For graphs this amounts to adding signs in the relations \S\ref{relpar}.
In particular, the following quadratic relations become anti--commutative:
\begin{eqnarray}
\scirct \ccirc{s'}{t'} &=&- \ccirc{s'}{t'} \scirct\\
\scirct\circ_{s't'}&=&-\circ_{s't'}\scirct\\
 \scirct \circ_{s't'}&=&-\ccirc{s'}{t'} \circ_{st}
\end{eqnarray}
Since \eqref{triangleeq} is not quadratic and hence the degree of a merger must be $0$ and
the relation does not get a sign
\begin{equation}
\scirct = \circ_{st}\mge{v}{w}
\end{equation}
Consequently, the following quadratic relations also remain without sign
\begin{eqnarray}
\mge{v}{w}\mge{v'}{w'}&=&\mge{v'}{w'}\mge{v}{w}\\
\scirct\mge{v}{w}& = &\mge{v}{w}\scirct\\
\circ_{st}\mge{v}{w}&=&\mge{v}{w}\circ_{st}
\end{eqnarray}
Isomorphisms also naturally have degree $0$ and hence there is no change in the relevant relation:
\begin{equation}
\phi\circ \sigma =-\sigma'\circ \phi'
\end{equation}

\subsubsection{Orders and Orientations}
In order to pictorially represent this, one can add decorations. This is very similar to the construction of ordered and oriented simplices, see e.g.\ \cite{munkres}. The first step is to give an order on all the edges of the ghost graph. The second step is to define orientations as orbits under even permutations. Finally one can impose the relation that two opposite orientations differ by a sign. Algebraically, one also uses the determinant line on the edges \cite{GKmodular}.
 It is only at this last step that the enrichment is needed. Furthermore  one can push this last step into the functor, that is only regard functors to Abelian $\C$ that take different change of orientations to sign changes.  These constructions are discussed in detail in \cite{feynman}.

\subsubsection{Graph Examples}

A list of examples
\begin{table}[!h]
\begin{tabular}{llll}
$\FF$&Feynman category for&condition on graphs +  additional decoration&\\
\hline
$\CCyclic^{odd}$&odd cyclic operads &trees + orientation of set of edges \\
$\modular^{odd}$&$\K$--modular&connected + orientation on set of edges \\
&&+ genus marking&\\
$\modular^{nc,odd}$&nc $\K$-modular& orientation on set of edges \\
&&+ genus marking&\\
$\dioperads^{\circlearrowleft odd}$&odd wheeled dioperads&directed graphs w/o parallel edges \\
&&+ orientations of edges&\\
$\props^{\circlearrowleft,ctd, odd}$&odd wheeled properads&connected directed graphs w/o parallel edges \\
&&+ orientation of set of edges&\\
$\props^{\circlearrowleft, odd}$&odd wheeled props &directed graphs w/o parallel edges \\
&&+ orientation of set of edges&\\
\end{tabular}
\caption{\label{oddtable}List of Feynman categories with conditions and decorations on the graphs}
\end{table}

\subsubsection{Suspension vs.\ odd}
In operad--lingo, one can suspend operads, etc.. On the Feynman category side this corresponds to certain twists.
I.e.\ there is a  twist $\Sigma$ and a $\Sigma$ twisted Feynman category $\FF_\Sigma$ such that  $\O\in \fopcatc$ iff the suspension $\Sigma\O\in\FF_{\Sigma}$--$\opcat_\C$. For general twistings of this type see \S\ref{twistsec}. These are equivalent to the odd version {\em if}  we are in the directed  case and there is a bijection between vertices and out flags,  see \cite{KWZ}. Even in the directed case, as explained in \cite{KWZ}.  the odd versions are actually more natural and yield the correct degrees in the Hochschild complex and correct signs and Master Equations, see \S\ref{univmasterpar}  below. A well known example for unexpected, but correct, signs is the Gerstenhaber bracket. It is {\em odd} Poisson.

In the same vein for the bar/cobar and Feynman transforms, it is not the suspended structures that are pertinent, but the odd structures, see  \S\ref{univmasterpar}.

\subsubsection{Examples}
\begin{enumerate}
\item
 Operads are very special, in the respect that their Feynman category is   equivalent to the one for their their odd version.
 \item The odd cyclic operads are equivalent to anti--cyclic operads.
\item For modular operads the suspended version is not equivalent to the odd versions a.k.a.\  $\K$--modular operads.
The difference is given by the twist $H_1(\gh(\phi))$.
\end{enumerate}

\subsection{Enriched versions. Details}

We can consider Feynman categories and target categories enriched over another monoidal category, such as $\mathcal{T}op$, $\Ab$ or $\dgVect$.
Note that  there are two cases. Either the enrichment is Cartesian, then we simply have to replace the free (symmetric) monidal category by the enriched version.
There is also a more categorical version of the definition with a condition going back to \cite{Getzler}. For that definition one simply replaces
all limits by indexed limits.
Or, the enrichment is not Cartesian, then we will replace the groupoid condition by an indexing just like above.

\subsubsection{Cartesian case: Categorical version} In \cite{feynman}
we proved that in the non--enriched case we can equivalently replace  (ii) by (ii').
\begin{itemize}
\item[(ii')] The pull-back of presheaves $\imath^{\otimes \wedge}\colon [\F^{op},Set]\to [\V^{\otimes op},Set]$
{\em restricted to representable presheaves} is monoidal.
\end{itemize}
This  then yields a definition in the Cartesian case if one replaces (iii) by the appropriate indexed limit condition.

\subsubsection{Non--Cartesian case indexed enrichment}
In the non--Cartesian case, the notion of groupoid ceases to make sense. The first option is to drop the groupoid condition and simply ask that the inclusion $\imath^{\otimes}$ is essentially surjective. This is possible and called a weak Feynman category, which is very close to the notion of  a pattern and explains that notion in more down to earth terms.
This is, however, not adequate for the bar/cobar and Feynman transforms or the twists.

The better notion is that of a Feynman category enriched over $\CalE$, indexed over another Feynman category $\FF$. The idea is that the Feynman category $\FF_\O$ for algebras over an operad $\O$ is a Feynman category enriched over $\C$ indexed over $\Surj$. The precise definition goes via enrichment functors, which are 2--functors.


In general, we will call the enrichment category $\CalE$.
 This is a monoidal category and hence can be thought of as a 2--category with one object, which we denote by $\underline{\CalE}$.
Here the 1-morphisms of $\underline{\CalE}$ are the objects of $\CalE$ with the composition being $\otimes$, the monoidal structure of $\CalE$. The 2--morphisms are then the 2--morphisms of $\CalE$, their horizontal composition being $\otimes$ and their vertical composition being $\circ$.
 Also, we can consider any category $\F$ to be a 2--category with the two morphisms generated by triangles of composable morphisms.

\begin{df}
Let $\FF$ be a Feynman category. An enrichment functor is a lax 2--functor $\D:\F\to \underline{\CalE}$ with the following properties
\begin{enumerate}
\item $\D$ is strict on compositions with isomorphisms.
\item $\D(\sigma)=\unit_{\CalE}$ for any isomorphism.
\item $\D$ is monoidal, that is $\D(\phi \otimes_{\F} \psi)=\D(\phi)\otimes_{\CalE}\D(\psi)$
\end{enumerate}
\end{df}

Given a monoidal category $\F$ considered as a 2--category and lax 2--functor $\D$ to $\underline{\CalE}$ as above, we define an enriched monoidal category $\F_{\D}$ as follows. The objects of $\F_{\D}$ are those of $\F$. The morphisms are given as

\begin{equation}
\label{edgedecoeq}
Hom_{\F_{\CalD}}(X,Y):=\bigoplus_{\phi\in Hom_{\F}(X,Y)}\D(\phi)
\end{equation}
The composition is given by
\begin{multline}
Hom_{\F_{\CalD}}(X,Y)\otimes Hom_{\F_{\CalD}}(Y,Z) =
\bigoplus_{\phi\in Hom_{\F}(X,Y)}\D(\phi)\otimes \bigoplus_{\psi\in Hom_{\F}(Y,Z)} \D(\psi)\\
\simeq
\bigoplus_{(\phi,\psi)\in Hom_{\F}(X,Y)\times Hom_{\F}(Y,Z)} \D(\phi)\otimes \D(\psi)
\stackrel{\bigoplus \D(\circ)}{\longrightarrow} \bigoplus_{\chi\in Hom_{\F}(X,Z)}\D(\chi)=Hom_{\F_{\D}}(X,Z)
\end{multline}
The image lies in the components $\chi=\psi\circ\phi$. Using this construction on $\V$, pulling back $\D$ via $\imath$, we obtain $\V_{\D}=\V_{\CalE}$, the freely enriched $\V$. The functor $\imath$ then is naturally upgraded to an enriched functor $\imath_{\CalE}:\V_{\D}\to \F_{\D}$.

\begin{df}
\label{enrichedfeydef}
Let $\FF$ be a Feynman category and let $\D$ be an enrichment functor.
We call $\FF_{\D}:=(\V_{\CalE},\F_{\D},\imath_{\CalE})$ a
Feynman category enriched
over $\CalE$ indexed by $\D$.
\end{df}

\begin{thm}
\label{enrichedtriplethm}
$\FF_{\D}$ is a weak Feynman category.  The forgetful functor from $\F_{\D}$-$\opcat$ to $\V_{\CalE}$-$\smodcat$ has a left adjoint and more generally push-forwards among indexed enriched Feynman categories exist. Finally there is an equivalence of
categories between algebras over the triple (aka.\ monad) $GF$ and $\F_{\D}$--$\opcat$.
\end{thm}

\begin{ex}
The freely enriched Feynman category. The functor $\D$ is simply the identity. This is the  triple
 $\FF_{\CalE}:=(\V_{\CalE},\F_{\CalE},\imath_{\CalE})$ where $\FF=(\V,\F,\imath)$ is a Feynman category and the subscript $\CalE$ means free enrichment.
\end{ex}

\begin{thm} The indexed enriched (over $\CalE$) Feynman category structures on a given FC $\FF$
are in 1--1 correspondence with  $\FF^{hyp}$-$\opcat$ and these are in 1--1 correspondence with enrichment functors.
\end{thm}
\begin{ex}
{\sc Twisted (modular) operads.}
Looking at $\FF=\modular$, we recover the notion of twisted modular operad. There is a twist for each hyper--operad $\D$.
We have the Feynman category $\modular_\D$. The triple then corresponds to $\mathbb{M}_\D$ in the notation of \cite{GKmodular}.
What we add is the descriptions (a) and (c), that is via compositions along graphs and generators and relations.
Here the graphs are actually decorated on the set of edges according to \eqref{edgedecoeq}. To see this one decomposes $\phi$ into simple edge or loop contractions as defined in \S\ref{zoopar}.
\end{ex}

\begin{ex}
Algebras over operads. In this case $\FF=\Surj$ and $\FF^{hyp,rd}=\operads_0$. An operad $\O\in \operads_0$-$\opcat_\C$ then gives an enrichment functor $\D_\O$ of $\Surj$. In particular $D_\O(n\twoheadrightarrow1)=\O(n)$ as in \S\ref{algebrasec}.
\end{ex}

\subsubsection{Coboundaries and $\V$--twists}
\label{twistsec}
 Coboundaries in the sense of \cite{GKmodular} are generalized to $\V$--twists.  Let $\fL\colon\V\to Pic(\CalE)$, that is the full subcategory of $\otimes$-invertible elements of $\CalE$.  A twist of a Feynman category indexed by $\D$ by $\fL$ is given by setting the new twist-system to be $\D_{\fL}(\phi)=
\fL(t(\phi))^{-1}\otimes \D(\phi)\otimes \fL(s(\phi))$.

The suspension  functor $s$ is such a coboundary twist, see \cite{GKmodular,KWZ}.
Here $\fL=s$ with $s(\ast_{(n-1)_+})=\Sigma^{2-n} sign_n$ in dg $\Vect$ for cyclic operads, or $s(\ast_{n_+})=\Sigma^{1-n} sign_n$ for operads, or in general $s (\ast_{(n-1)_+})=\Sigma^{-2(g-1)+n}sign_n$ where $\Sigma$ is the suspension and $sign_n$ is the sign representation, see \cite{KWZ} for a detailed explanation.

\subsubsection{Odd versions and shifts}
\label{oddtwistpar}

Given a well-behaved presentation of a Feynman category (generators+relations for the morphisms) we can define an odd version which is enriched over $\Ab$ by giving a twist.  To obtain the odd versions, we use $\D(\phi)=\det(Edges(\gh(\phi))$.
In the cyclic case, an example are anti--cyclic operads and the theory of modular operads this twist is called $\K$. It is {\em not} a coboundary in general. Rather up to the suspension coboundary and the shift coboundary, this twist is a twist by $H_1(\gh)$ in the modular case, see \cite{GKmodular,KWZ} for details.

\subsection{Feynman  level category $\FF^+$, hyper category $\FF^{hyp}$ and its reduction $\FF^{hyp,rd}$.}
\subsubsection{Feynman  level category $\FF^+$}
\label{Fplussec}
Given a Feynman category $\FF$, and a choice of basis for it, we will define its
Feynman level  category  $\FF^+=(\V^+,\F^+,\imath^+)$ as follows.
The underlying objects of $\F^+$ are the morphisms of $\F$. The morphisms
of $\F^+$ are given as follows: given $\phi$ and $\psi$, consider their decompositions
\begin{equation}
\xymatrix
{
X \ar[rr]^{\phi}\ar[d]_{\sigma}^{\simeq}&& Y\ar[d]^{\hat\sigma}_{\simeq} \\
 \bigotimes_{v\in I}\bigotimes_{w\in I_v} \ast_w\ar[rr]^{\bigotimes_{v\in I}\phi_{v}}&&\bigotimes_{v\in I} \ast_v
}
\quad
\xymatrix
{
X' \ar[rr]^{\psi}\ar[d]_{\tau}^{\simeq}&& Y'\ar[d]^{\hat\tau}_{\simeq} \\
 \bigotimes_{v'\in I'}\bigotimes_{w'\in I'_{v'}} \ast_{w'}\ar[rr]^{\bigotimes_{v'\in I'}\psi_{v'}}&&\bigotimes_{v'\in I'} \ast_{v'}.
}
\end{equation}
where we have dropped the $\imath$ from the notation, $\sigma,\hat\sigma,\tau$ and $\hat\tau$ are given  by the choice of basis
and the partition $I_v$ of the index set for $X$ and $I'_{v'}$ for the index set of $Y$
is given by the decomposition of the morphism.

 A morphism from $\phi$ to $\psi$ is a two level partition of $I: (I_{v'})_{v'\in I'}$, and partitions of $I_{v'}:(I^1_{v'}\dots,I_{v'}^{k_{v'}})$ such that
 if we set $\phi_{v'}^i:=\bigotimes_{v\in I_{v'}^i}\phi_v$ then $\psi_{v'}=\phi_{v'}^k\circ\dots\circ \phi_{v'}^1$.

 To compose two morphisms $f\colon\phi\to \psi$ and $g\colon\psi\to \chi$, given by partitions of $I:(I_{v'})_{v'\in I'}$ and of the $I_{v'}:(I^1_{v'}\dots,I_{v'}^{k_{v'}})$ respectively of
 $I':(I'_{v''})_{v''\in I''}$ and the $I_{v''}:(I^{\prime 1}_{v''}\dots,I_{v''}^{\prime k_{v''}})$,
 where $I''$ is the index set in the decomposition of $\chi$, we set the compositions to be
 the partitions of $I:(I_{v''})_{v''\in I''}$ where  $I_{v''}$ is the set partitioned by
 $(I_{v'})_{v'\in I^{\prime j}_{v''},j=1,\dots, k_{v''}}$.
 That is, we replace each morphism $\psi_{v'}$ by the chain $\phi_1^{v'}\circ\dots\circ\phi_k^{v'}$.

\begin{figure}[h]
    \centering
   \includegraphics[width=0.8\textwidth]{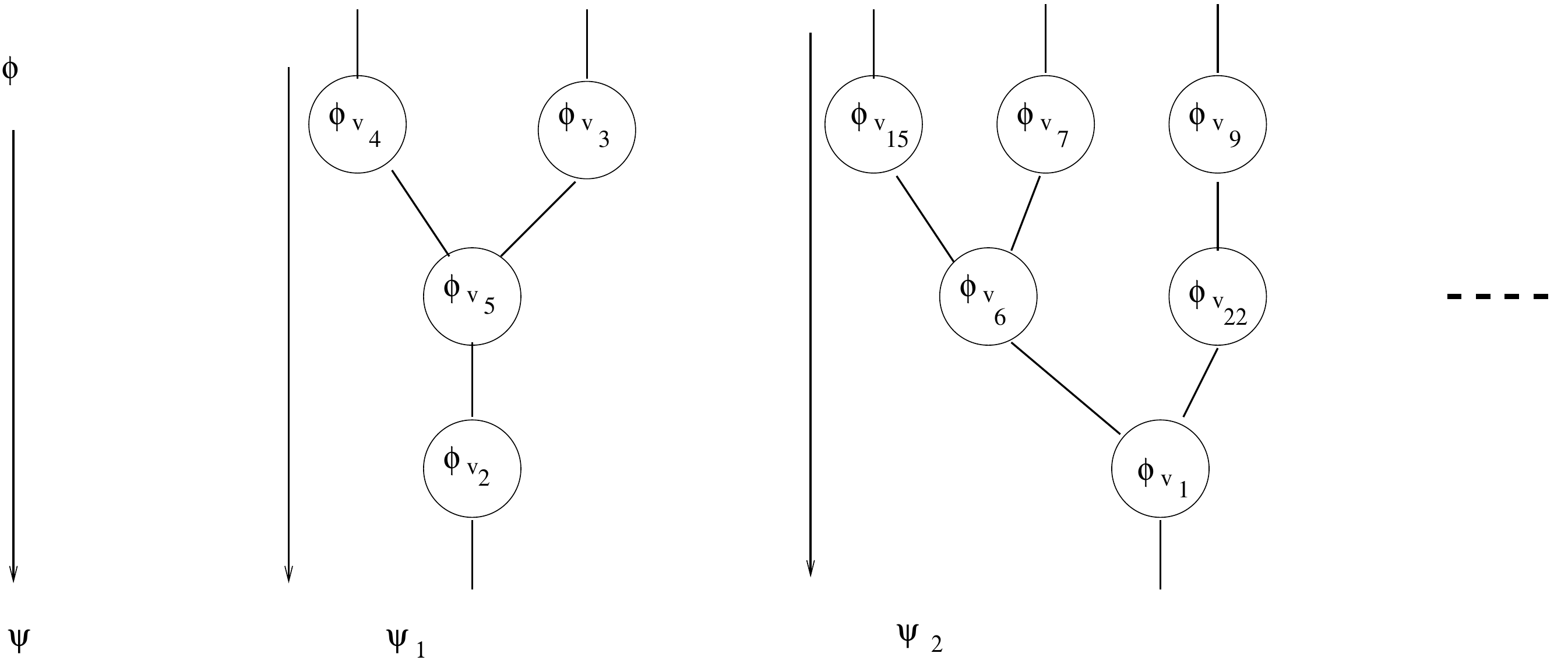}
    \caption{The level forest picture for morphisms in $\FF^+$. Indicated is a morphism from $\phi\simeq\bigotimes_v\phi_v$ to $\Psi\simeq \bigotimes_{i}\Psi_i$}
    \label{Fplusfig}
\end{figure}

Morphisms alternatively correspond to  rooted forests of level trees thought of as flow charts, see Figure \ref{Fplusfig}. Here the vertices are decorated by the $\phi_v$ and the composition along the rooted  forest is $\psi$.
There is exactly one tree $\tau_{v'}$ per $v'\in I'$ in the forest and accordingly
the composition along that tree is $\psi_v'$.

Technically,   the vertices are the $v\in I$.
The flags are the union $\amalg_v \amalg_{w\in I_v} \ast_w\amalg \amalg_{v\in I}\ast_v$
with the value of $\del$ on $\ast_w$ being $v$ if  $w\in I_v$ and
$v$ on $\ast_v$ for $v\in I$. The orientation at each vertex is given by the target being out. The involution $\imath$  is given by matching source and target objects of the various $\phi_v$.
The level structure of each tree is given by the partition  $I_{v'}$.
 The composition is the composition of rooted  trees by gluing trees  at all vertices ---that is we blow up the vertex marked by $\psi_{v'}$ into the tree $\tau_{v'}$.

\subsubsection{$\F^+$-$\opcat$.}
After passing to the equivalent strict Feynman category, an element $\CalD$ in $\F^+$-$\opcat$ is a symmetric monoidal functor that has values on each morphism $\CalD(\phi)=\bigotimes\CalD(\phi_v)$ and has composition maps $\CalD(\phi_0\otimes \phi)\to \CalD(\phi_1)$ for each decomposition $\phi_1=\phi\circ\phi_0$. Further decomposing $\phi=\bigotimes \phi_v$ where the decomposition is according to the target of $\phi_0$, we obtain morphisms
\begin{equation}
\label{hypercompeq}
 \CalD(\phi_0)\otimes \bigotimes_v\CalD(\phi_v)\to \CalD(\phi_1)
 \end{equation}

It is enough to specify these functors for $\phi_1\in (\F \downarrow \V)$ and then check associativity for triples.

\begin{ex}
If we start from the tautological  Feynman category on the trivial category $\FF=(1,1^{\otimes}, \imath)$ then $\FF^+$ is the Feynman category $\FSurj$ of surjections. Indeed the possible trees are all linear, that is only have 2--valent vertices, and there is only one decoration. Such a  rooted tree is specified by its total length $n$ and the permutation which gives the bijection of its vertices with the set $n_i$. Looking at a forest of these trees we see that we have the natural numbers as objects with morphisms being surjections.
\end{ex}
\begin{ex}
We also have $\FSurj^+=\operads_{May}$, which is the Feynman category for  May operads. Indeed the basic maps (\ref{hypercompeq}) are precisely the composition maps $\gamma$. To be precise, these are May operads without units.

\end{ex}

\subsubsection{Feynman hyper category $\FF^{hyp}$}
\label{hypersec}
There is a ``reduced'' version of $\FF^+$ which is central to our theory of enrichment. This is the universal Feynman category through which any functor $\CalD$ factors if it satisfies the following restriction:
{\it  $\CalD(\sigma)\simeq\unit$ for any isomorphism $\sigma$ where $\unit$ is the unit of the target category $\CalC$.}

For this, we invert the morphisms corresponding to composing with isomorphisms, see \cite{feynman} for details.

\subsubsection{$\F^{hyp}$-$\opcat$}
\label{hypopsec}
An element $\D\in \F^{hyp}$-$\opcat$ corresponds to the data of functors from $Iso(\F\downarrow\F)\to \C$ together with morphisms (\ref{hypercompeq}) which are associative and satisfy the condition that all the following diagrams commutes:
\begin{equation}
\label{isoactioneq}
\xymatrix{
\CalD(\phi)\otimes \bigotimes_v \CalD(\sigma_v)\ar[rr]^{\CalD(\tau)}_{\sim}&&\CalD({}^ {\boldsymbol \sigma}\phi)\\
\CalD( \phi)\otimes \bigotimes_v \unit  \ar[u]^{\sim}&& \ar[ll]_{id \otimes \bigotimes_v r^{-1}_{\unit}}^{\sim}\D(\phi)\ar[u]_{\CalD( \boldsymbol \sigma)}
}
\end{equation}
see \cite{feynman} for details.

\begin{ex}
The paradigmatic examples are hyper--operads in the sense of \cite{GKmodular}. Here $\FF=\modular$ and $\FF^{hyp}$  is the Feynman category for hyper--operads.
\end{ex}

\subsubsection{A reduced version $\FF^{hyp,rd}$}  One may define $\FF^{hyp,rd}$, a Feynman subcategory of $\FF^{hyp}$ which is equivalent to it by letting $\F^{hyp, rd}$ and $\V^{hyp, rd}$ be the respective subcategories whose objects are morphisms that do not contains isomorphisms in their decomposition. In view of the isomorphisms $\emptyset\to \sigma$ this is clearly an equivalent subcategory.  In particular, the respective categories of $\opcat$ and $\op{M}ods$ are equivalent.

The  morphisms are  described by rooted forests of  trees whose vertices are decorated by the $\phi_v$ as above --none of which is an isomorphism--, with the additional decoration of an isomorphism
per edge and tail. Alternatively, one can think of the decoration as a black 2-valent vertex.
Indeed, using maps from $\emptyset \to \sigma$, we can introduce as many isomorphisms as we wish. These give rise to 2--valent vertices, which we mark black. All other vertices remain labeled by $\phi_v$.  If there are sequences of such black vertices, the corresponding morphism is isomorphic to the morphism resulting from composing the given sequence of these isomorphisms.

\begin{ex}
For $\FF^{hyp,rd}_{surj}=
\operads_0$, the Feynman category whose morphisms are trees with at least trivalent vertices (or identities) and whose $\opcat$ are operads whose $\O(1)=\unit$. Indeed the basic non--isomorphism
 morphisms are the surjections $\underline{n}\to \underline{1}$, which we can think of as rooted corollas. Since for any two singleton sets there is a unique isomorphism between them, we can suppress the black vertices in the edges. The remaining information is that of the tails, which is exactly the map $\phi^F$ in the morphism of graphs.
\end{ex}
\begin{ex}
For the trivial Feynman category, we obtain back the trivial Feynman category as the reduced hyper category, since the trees all collapse to a tree with one black vertex.
\end{ex}

\subsection{Free monoidal construction $\F^{\boxtimes}$}
\label{freemsec}
Sometimes it is convenient to construct a new Feynman category from a given one whose vertices
are the objects of $\F$.
Formally, we set $\FF^{\boxtimes}=(\V^{\otimes},\F^{\boxtimes},\imath^{\otimes})$ where $\F^{\boxtimes}$ is the free monoidal
category on $\F$ and we denote the ``outer'' free monoidal structure by $\boxtimes$.
This is again a Feynman category.
There is a functor $\mu:\F^{\boxtimes}\to \F$ which sends $\boxtimes_i X_i\mapsto \bigotimes_i X_i$
and by definition $Hom_{\F^{\boxtimes}}({\bf X}=\boxtimes_i X_i,{\bf Y}=\boxtimes_iY_i)=\bigotimes_{i}Hom_{\F}(X_i,Y_i)$. The only way that the index sets can differ, without the Hom--sets being empty, is if some of the factors are $\unit\in \F^{\boxtimes}$. Thus the one--comma generators are simply the elements of $Hom_\F(X,Y)$.
Using this identification one obtains: $Iso(\F^{\boxtimes})\simeq Iso(\F)^{\boxtimes}\simeq( \V^{\otimes})^{\boxtimes}$.
The factorization and size axiom follow readily from this description.

\begin{prop}
$\F^{\boxtimes}$-$\opcat_{\C}$ is equivalent
to the category of functors (not necessarily monoidal)
$Fun(\F,\C)$.
\end{prop}

\begin{ex}
Examples are  $FI$ modules and (crossed) simplicial objects for the free monoidal Feynman categories for $FI$ and
$\mathbf \Delta_+$ where for the latter one uses the non--symmetric version.
\end{ex}

\subsection{NC--construction}
\label{ncsec}
For any Feynman category one can define its nc (non--connected) version.
It plays a crucial role in physics and mathematics and manifests itself through the BV equation \cite{KWZ}. Namely, for the operator $\Delta$ in the case of modular operads to become a differential, one needs a multiplication. This, on the graph level, is given by disjoint union for the one--comma generators. This amounts to dropping the condition of connectedness.  Astonishingly this works in full generality for any Feynman category.

Let $\FF=(\V,\F,\imath)$, then we set $\FF^{nc}=(\V^{\otimes},\F^{nc},\imath^{\otimes})$ where $\F^{nc}$ has objects $\F^{\boxtimes}$, the free monoidal product. We however add more morphisms. The one--comma generators will be
$Hom_{\F^{nc}}({\bf X},Y):=Hom_\F(\mu({\bf X}),Y)$,
where for ${\bf X}=\boxtimes_{i\in I}X_i$, $\mu({\bf X})=\bigotimes_{i\in I}X_i$.
This means that for ${\bf Y}=\boxtimes_{j\in J}Y_j$,
$Hom_{\F}({\bf X},{\bf Y})\subset Hom_{\F}(\mu({\bf X}),\mu({\bf Y}))$, includes only those morphisms for which there is a partition $I_j, j\in J$ of $I$ such that the morphism factors through $\bigotimes_{j\in J} Z_j$ where $Z_j\stackrel {\sigma_j}{\to}\bigotimes_{k\in I_j} X_k$ is an isomorphism. That is $\psi=\bigotimes_{j\in J} \phi_j\circ \sigma_j$ with $\phi_j:Z_j\to Y_j$. Notice that there is a map of ``disjoint union'' or ``exterior multiplication''  given by $\mu:X_1\boxtimes X_2\to X_1\otimes X_2$ via $id\otimes id$.

\begin{ex}
The terminology ``non--connected'' has its origin in the graph examples.
Examples can be found in \cite{KWZ}, where also a box--picture for graphs is presented. The connection is that morphisms in $\F^{nc}$ have an underlying graph that is disconnected and the connected components are those of the underlying $\F$.

\end{ex}

\begin{prop}\cite{feynman}
There is an equivalence of categories between $\F^{nc}$-$\opcat_\C$ and
symmetric lax monoidal functors $Fun_{lax\;\otimes}(\F,\C)$.
\end{prop}

Using lax--monoidal functors, is also a way to deal with algebras over operads whose $\O(1)$ has isomorphisms.

%
%
%
%
%
%

\section{Universal operations, Transforms and Master equations}
\label{univmasterpar}
\subsection{Universal operations}
\subsubsection{Universal operations for Operads, etc.}
A well known result in operad theory is that for an operad $\O$ there is an odd Lie bracket defined on $\bigoplus \O(n)$ \cite{G}.
This actually descends to coinvariants  $\bigoplus \O(n)_{\SS_n}$ \cite{KapMan}. For anti--cyclic operads there is again an odd Lie bracket on the
coinvariants $\bigoplus \O((n))_{\SS_n}$ with lifts to the smaller coinvariants w.r.t.\ the cyclic groups $C_n$, namely on $\bigoplus_n \O((n))_{C_n}$ \cite{KWZ}.
Similarly there are operations $\Delta$ on $\bigoplus \O((n,g))_{\SS_n}$ for modular operads \cite{KWZ,Bar}. Here we show that these operations can be understood purely from the Feynman category and we can explain why exactly these operations turn up in the master equations.

\subsubsection{Cocompletion}
Let $\hat \F$ be the cocompletion of $\F$. This is monoidal with the monoidal structure given by the Day convolution $\day$. If $\C$ is cocomplete then $\O\in \opcat$ factors:
$$\xymatrix{
\F\ar[dr]^{\jmath}\ar^{\O}[rr]&&\CalC\\
&\hat \F\ar^{\hat \O}[ur]&\\
}$$

\begin{thm} Let $\Eins:=\colim_{\V}\jmath \circ \imath \in \hat \F$ and let $\F_{\V}$ the symmetric monoidal subcategory generated by $\Eins$.
Then $\FF_{\V}:=(\F_\V,\Eins,\imath_{\V})$ is a Feynman category. (This gives an underlying operad of universal operations).
\end{thm}

If $\CalE$ is Abelian,
we say $\FFV$ is weakly generated by morphisms $\phi\in \Phi$ if the summands of the components $[\phi_{X_{\bf j},i}]$
generate the morphisms of $\FFV$. Here different summands are indexed by different isomorphism classes of morphisms.

\subsubsection{Example: Operads}
$\operads$ the Feynman category for operads, $\C=\dgVect$.

Then $\hat \O(\Eins)=\bigoplus_n\O(n)_{\SS_n}$ and the Feynman category is (weakly) generated by
$\circ:=[\sum \circ_i]$. (This is a two--line calculation).
This gives rise to the Lie bracket by using the anti--commutator.
It lifts to the non-Sigma case  along the forgetful  $\operads^{\neg\Sigma}\to \operads$ and gives the pre--Lie structure on $\bigoplus_n\O(n)$, which
goes back to \cite{G}. In  \cite{KapMan} it was shown that the pre--Lie structure descends to the coinvariants. In \cite{KWZ} it is argued that the pre--Lie structure lives naturally on the coinvariants and lifts to the invariants.

In general these kinds of lifts are possible if there is a non--Sigma version.

\subsubsection{Example: Odd/anti--cyclic Operad}
The universal operations are (weakly) generated by a Lie bracket. $[ \, ,\, ]:=[\sum_{st}\circ_{st}]$, (see \cite{KWZ}).
This actually lifts to cyclic coinvariants (non--sigma cyclic operads) that is along the map $\CCyclic^{odd,pl}\to \CCyclic^{odd}$.
Here we also see that one cannot expect a further lift, since the planar version for $\CCyclic^{odd}$ still has a non discrete $\V$.

\subsubsection{The three geometries of Kontsevich}
The endomorphism operad $End(V)$ for a symplectic vector space is anti--cyclic.
Any tensor product: $(\O \otimes \P)(n):=\O(n)\otimes \P(n)$ with $\O$ a cyclic operad and $\P$ an  anti--cyclic operad is anti--cyclic and hence has
the odd Lie bracket discussed above.

Fix $V^n$ $n$--dim symplectic $V^{n}\to V^{n+1}$.
For each $n$ get Lie algebras
\begin{enumerate}
\item $Comm\otimes End(V^n)$
\item  $Lie\otimes End(V^n)$
\item $Assoc \otimes End(V^n)$
\end{enumerate}

Taking the limit as $n\to \infty$ one obtains the formal geometries of \cite{kontsevichthree,ConantVogtmann}.

Our construction is more general and works for any anti--cyclic operad.
For instance another family of Lie algebras can be obtained as follows, \cite{KWZ}.
Let $V^n$ be a vector space with a symmetric non-degenerate form. $End(V)$ is a cyclic operad. Since the $PreLie$ operad is anti--cyclic \cite{cha}, for each $n$ we get a Lie algebra $PreLie\otimes End(V)$. It is not known what geometry we get when we take the limit as $n\to \infty$.

\subsubsection{Further examples}
For further examples, see Table \ref{univtable}.
\begin{table}[h]
\begin{tabular}{lllll}
$\FF$&Feynman cat for&$\FF,\FF_{\V}$,$\FF_{\V}^{nt}$&weak gen.~subcat.&\\
\hline
$\operads$&Operads&rooted trees&$\Fprelie$\\
$\operads^{odd}$&odd operads&rooted trees + orientation & odd pre-Lie\\
 &&of set of edges\\
$\operads^{\neg\Sigma}$&non-Sigma operads &planar rooted trees & all $\circ_i$ operations\\
$\operads_{mult}$&Operads with mult.&b/w rooted trees&pre-Lie + mult.\\
$\CCyclic$&cyclic operads&trees& com. mult.&\\
$\CCyclic^{odd}$&odd cyclic operads &trees + orientation & odd Lie&\\
&&of set of edges\\
$\modular^{odd}$&$\K$--modular&connected + orientation & odd dg Lie \\
&& on set of edges \\
$\modular^{nc,odd}$&nc $\K$-modular& orientation on set of edges & BV\\
$\dioperads$&Dioperads&connected directed graphs w/o&Lie--admissible\\
&& directed  loops or parallel edges&\\
\end{tabular}
\caption{\label{univtable}Here $\FF_{\V}$ and $\FF_{\V}^{nt}$ are
given as $\FF_{\O}$ for the operad $\O$, the composition as discussed being insertion. The former is for the type of graph with unlabelled tails and
the latter for the version with no tails.}
\end{table}

\subsection{Transforms \& Master Equations}
There are three transforms we will consider: the bar--, the cobar transform and the Feynman transform aka.\ dual transform.

\subsubsection{Motivating example: Algebras}
If $A$ is an associative algebra, then the bar transform is the dg--coalgebra given  by the free  coalgebra $BA=T\Sigma^{-1}\bar A$ together with co--differential from algebra structure. The usual notation for an element in $BA$ is $a_0|a_1|\dots|a_n$.

Likewise let $C$ be an associative co--algebra.  The co--bar transform is the dg--algebra $\Omega C:=Free_{alg}(\Sigma^{-1}\bar C)$ together with a differential coming from co--algebra structure.
The bar--cobar transform
$\Omega BA$ is a  resolution of $A$.

For the Feynman transform consider $A$ a finite--dimensional algebra or graded algebra with finite dimensional pieces and let $\check A$ be its dual co--algebra.
Then the dual or Feynman transform of $A$ is $FA:=\Omega \check A$ + differential from multiplication.  Now,  the double Feynman transform $FFA$ a resolution.

\subsubsection{Transforms}
\begin{figure}
    \centering
    \includegraphics[width=0.5\textwidth]{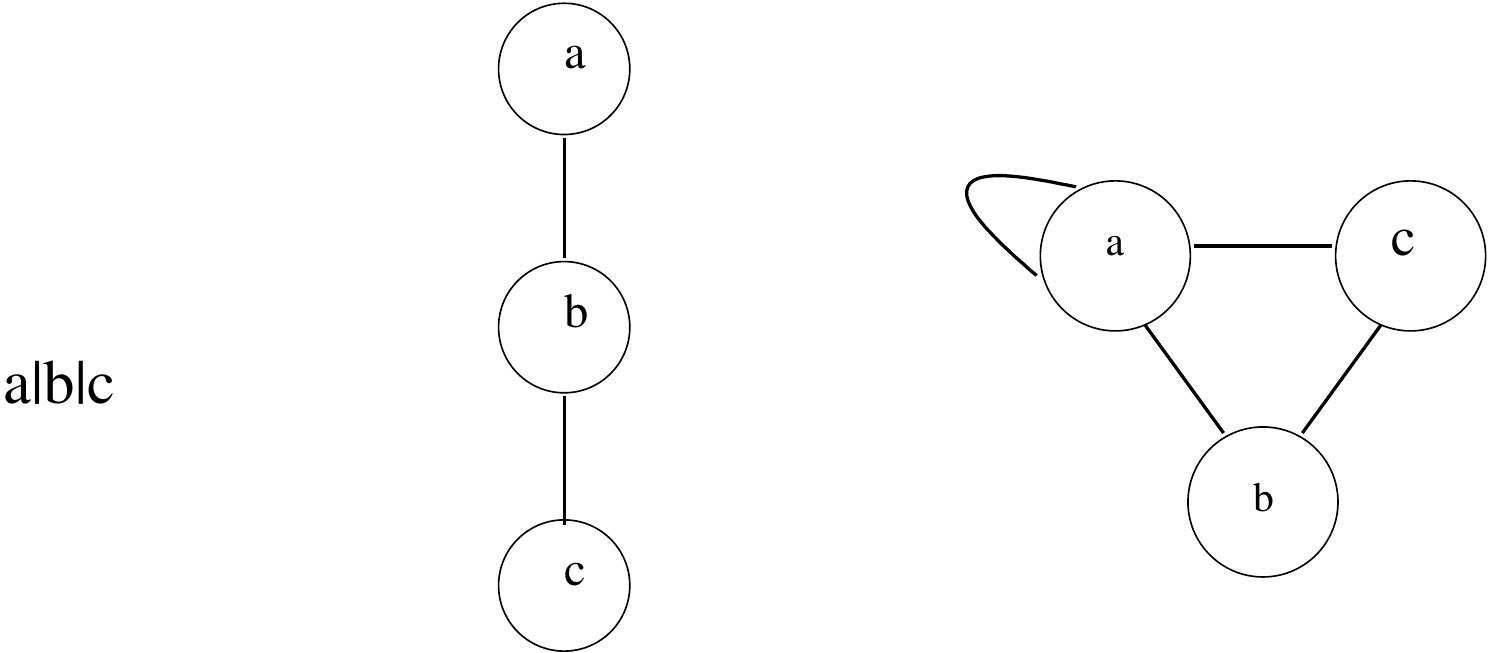}
    \caption{The sign mnemonics for the bar construction, traditional version with the symbols $|$ of degree $1$, the equivalent linear tree with edges of degree $1$,  and a more general graph with edges of degree $1$. Notice that in the linear case there is a natural order of edges, this ceases to be the case for more general graphs}
    \label{barfig}
\end{figure}

These transforms take $\O\in \fopcatc$ and transform it to an $\oper$ for the odd version of the Feynman category $\FF^{odd}$ either in $\C^{op}$ or $\C$.
All these are free constructions, which, however, also have the extra structure of an additional (co)differential. Thus the resulting Feynman category is actually enriched over chain complexes and one can start out there as well. Furthermore, for the (co)differential to work, we have to have signs. These are exactly what is provided by the odd versions. In order to be able to define the transforms, one has to fix an odd version $\FF^{odd}$ of $\FF$, just as in \S\ref{oddgraphpar}. This is analogous to the suspension in the usual bar transforms. In fact, the following is more natural, see \cite{feynman, KWZ}. The degree is $1$ for each bar and in the graph case the edges get degree $1$; see   Figure \ref{barfig}.
We can generalize the construction of $\FF^{odd}$ to so--called well--presented Feynman categories, see below and \cite{feynman}.
In this case, we can define the  transformations for elements of $\opcat$.

The Feynman transform is of particular interest. Since the construction is free, any $\V\in \smodcat$ will yield an $\oper$.
On the other hand, this need not be compatible with the dg structure. It turns out that it is, if it satisfies a Master Equation.

The transforms are of interest in themselves, but one common application is that the bar-cobar transform as well as the double Feynman transform give a ``free'' resolution. In general, of course, ``free'' means co-fibrant. For this kind of statement one needs a Quillen model structure, which is provided in \S\ref{modelpar}.

\begin{rmk}
As before one can ask the question of how much of the structure of these transforms can be pulled back to the Feynman category side. The answer is: ``Pretty much all of it''. We shall not discuss this here, but it can be found in \cite{feynman}.
\end{rmk}

\subsubsection{Presentations } In order to define the transforms,  we have to give what is called an ordered presentation \cite{feynman}.
Rather then giving the technical conditions, we will consider the graph case  and show these structures in this case.

\subsubsection{Basic example $\GG$}
In $\GG$ the presentation comes from the following set of morphisms $\Phi$
\begin{enumerate}
\item There are 4 types of basic morphisms: Isomorphisms, simple edge contractions, simple loop contractions and mergers. Call this set $\Phi$.
\item
 These morphisms generate all one--comma generators upon iteration. Furthermore, isomorphisms act transitively on the other classes.  The relations on the generators are  given by commutative diagrams.
 \item
 The relations are  quadratic for edge contractions as are the relations involving isomorphisms.
 Finally there is a non--homogenous relation coming from  a simple merger and a loop contraction being equal to an edge contraction.

\item  We can therefore assign degrees as $0$ for isomorphisms and mergers, $1$ for edge or loop contractions and split $\Phi$ as $\Phi^0\amalg \Phi^1$. This gives a degree to any morphism.
\end{enumerate}

Up to isomorphism any morphism of degree $n$ can be written in $n!$ ways up to morphisms of degree $0$. These are the enumerations of the edges of the ghost graph.

There is also a standard order in which isomorphisms come before mergers which come before edge contractions as in \eqref{bmdecompeq}.
This gives an ordered presentation.

In general, an ordered presentation is a set of generators $\Phi$ and extra data such as the subsets $\Phi^0$ and $\Phi^1$; we refer to \cite{feynman} for details.

\subsubsection{Differential} Given a
$d_{\Phi^1}=\sum_{[\phi_1]\in \Phi^1/\sim} \phi_1\circ$ defines an endomorphism on the Abelian group generated by the isomorphism classes morphisms.
The non--defined terms are set to zero. $\Phi^1$ is called resolving if this is a differential.

In the graph case, this amounts to the fact that for any composition of edge contractions $\phi_e \circ \phi_{e'}$, there is precisely another pair of edge contractions
$\phi_{e''} \circ \phi_{e'''}$ which contracts the edges in the opposite order.

This differential will induce differentials for the transforms, which we call by the same name. We again refer to \cite{feynman} for details.

\subsubsection{Setup} $\FF$ be a Feynman category enriched over $\Ab$ and with an ordered presentation and let $\FF^{odd}$ be its corresponding odd version. Furthermore let $\Phi^1$ be a resolving subset of one-comma generators and let $\op{C}$ be an additive category, i.e.\ satisfying the analogous conditions above. In order to give the definition, we need a bit of preparation.
Since $\V$ is a groupoid, we have that $\V\simeq \V^{op}$. Thus, given a functor $\Phi:\V\to \C$, using the equivalence we get a functor from $\V^{op}$ to $\C$ which we denote by $\Phi^{op}$. Since the bar/cobar/Feynman transform adds a differential, the natural target category from $\fops$ is not $\C$, but complexes in $\C$, which we denote by $Kom(\C)$. Thus any $\O$ may have an internal differential $d_\O$.

\subsubsection
{ The bar construction} This is the functor
\begin{equation*}
\Bar \colon \fops_{Kom(\C)}\to \foddops_{Kom(\C^{op})}
\end{equation*}
\begin{equation*}
\Bar(\O):=\imath_{\FF^{odd} \; *}(\imath_{\FF}^*(\O))^{op}
\end{equation*}
together with the differential $d_{\op{O}^{op}}+d_{\Phi^1}$.

\subsubsection{  The cobar construction} This is the functor
\begin{equation*}
\Cobar \colon \foddops_{Kom(\C^{op})}\to \fops_{Kom(\C)}
\end{equation*}
\begin{equation*}
\Cobar(\op{O}):=\imath_{\FF \; *}(\imath^\ast_{\FF^{odd}}(\op{O}))^{op}
\end{equation*}
together with the co-differential $d_{\op{O}^{op}}+d_{\Phi^1}$.

\subsubsection{Feynman transform}
 Assume there is a duality equivalence $\vee\colon \CalC\to \CalC^{op}$.
The Feynman transform is a pair of functors, both denoted $\FT$,
\begin{equation*}
\FT\colon \fops_{Kom(\C)}  \leftrightarrows \foddops_{Kom(\C)}\colon \FT
\end{equation*}
defined by
\begin{equation*}
 \FT(\O):=\begin{cases} \vee\circ \Bar(\O) & \text{ if } \O \in \fops_{Kom(\C)} \\ \vee\circ \Cobar(\O) & \text{ if } \O \in \foddops_{Kom(\C)}
\end{cases}
\end{equation*}

\begin{prop} The bar and cobar construction form an adjunction.
\begin{equation*}
\adj{\Cobar}{\foddops_{Kom(\op{C}^{op})}}{\fops_{Kom(\op{C})}}{\Bar}
\end{equation*}
\end{prop}

The quadratic relations in the graph examples are a feature that can be generalized to the notion of {\em cubical} Feynman categories. The name reflects the fact that in the graph example the $n!$ ways to decompose a morphism whose ghost graph is connected and has $n$ edges into simple edge contractions correspond to the edge paths of $I^n$ going from $(0,\dots,0)$ to $(1,\dots, 1)$. Each  edge flip in the path represent one of the quadratic relations and furthermore the $\SS_n$ action on the coordinates is transitive on the paths, with transposition acting as edge flips.

This is a convenient generality in which to proceed.

\begin{thm} Let $\FF$ be a cubical Feynman category and $\op{O}\in \fops_{Kom(\C)}$.  Then the counit $\Cobar\Bar(\op{O})\to\op{O}$ of the above adjunction is a levelwise quasi-isomorphism.
\end{thm}

\begin{rmk}
In the case of $\mathcal{C}=dgVect$, the Feynman transform can be intertwined with the aforementioned push-forward and pull-back operations to produce new operations on the categories $\mathcal{F}-\mathcal{O}ps_\mathcal{C}$.  A lifting (up to homotopy) of these new operations to $\mathcal{C}=Vect$ is given in \cite{Ward17}.  In particular this result shows how the Feynman transform of a push-forward (resp. pull-back) may be calculated as the push-forward (resp. pull-back) of a Feynman Transform.  One could thus assert that the study of the Feynman transform belongs to the realm of Feynman categories as a whole and not just to the representations of a particular Feynman category.
\end{rmk}

\subsection{Master equations}
In \cite{KWZ}, we identified the common background  of master equations that had appeared throughout the literature for operad--like objects and extended them to all graphs examples. An even more extensive theorem for Feynman categories can also be given.

The Feynman transform is quasi--free. An algebra over $F\O$ is dg--if and only if it satisfies the relevant Master Equation. First, we have the tabular theorem from \cite{KWZ} for the usual suspects.

\begin{thm}\label{methm}(\cite{Bar},\cite{MerkVal},\cite{wheeledprops},\cite{KWZ})  Let $\op{O}\in \fopsc$ and $\op{P}\in \foddops_\op{C}$ for an $\F$ represented in Table $\ref{MEtable}$.  Then there is a bijective correspondence:
\begin{equation*}
Hom(\FT(\op{P}),\op{O})\cong ME(lim_\V(\op{P} \tensor\op{O}))
\end{equation*}


\end{thm}
Here ME is the set of solutions of the appropriate master equation set up in each instance.

\begin{table}[htb] \centering
\begin{tabular}{p{2.2cm}||l|l}
Name of  $\fopsc$& Algebraic Structure of $F\O$& Master Equation (ME)\\ \hline\hline
operad,\cite{GJ}& odd pre-Lie & $d(-)+ -\circ- =0$  \\ \hline
cyclic operad \cite{GKcyclic}  & odd Lie & $d(-)+ \frac{1}{2}[-,-] =0$  \\ \hline
modular operad \cite{GKmodular}& odd Lie + $\Delta$& $d(-)+ \frac{1}{2}[-,-]+\Delta(-) =0$  \\ \hline
properad  \cite{Vallette}& odd pre-Lie & $d(-)+ -\circ- =0$  \\ \hline
wheeled properad \cite{wheeledprops} & odd pre-Lie + $\Delta$ & $d(-)+ -\circ- +\Delta(-) =0$  \\ \hline
wheeled prop \cite{KWZ} & dgBV & $d(-)+ \frac{1}{2}[-,-] +\Delta(-) =0$  \\
\end{tabular}
\caption{
\label{MEtable} Collection of Master Equations for operad--type examples}
\end{table}

 With Feynman categories this tabular theorem can be compactly written and generalized.  The first step is the realization that the differential specifies a natural operation, in the above sense, for each arity $n$. Furthermore, in the Master Equation there is one term form each generator of $\Phi^1$ up to isomorphism.
 This is immediate from comparing Table \ref{MEtable} with Table \ref{univtable}.
  The natural operation which lives on a space associated to an $\op{Q}\in \fopcat$ is denoted $\Psi_{\op{Q},n}$ and is formally defined as follows:

\begin{df}
For a Feynman category $\FF$ admitting the Feynman transform and for $\op{Q}\in\fopsc$ we define the formal master equation of $\FF$ with respect to $\op{Q}$ to be the completed cochain $\Psi_{\op{Q}}:= \prod \Psi_{\op{Q},n}$.  If there is an $N$ such that $\Psi_{\op{Q},n}=0$ for $n>N$, then we define the master equation of $\FF$ with respect to $\op{Q}$ to be the finite sum:
\begin{equation*}
d_{\op{Q}}+\ds\sum_{n}\Psi_{\op{Q},n} = 0
\end{equation*}
We say $\alpha\in lim_\V(\op{Q})$ is a solution to the master equation if $d_\op{Q}(\alpha)+\sum_{n}\Psi_{\op{Q},n}(\alpha^{\tensor n}) = 0$, and we denote the set of such solutions as $ME(lim_\V(\op{Q}))$.
\end{df}
Here the first term is the internal differential and the term for $n=1$ is the differential corresponding to $d_{\Phi^1}$, where $\Phi^1$ is the subset of odd generators.

\begin{thm} Let $\op{O}\in \fopsc$ and $\op{P}\in \foddops_\op{C}$ for an $\F$ admitting a Feynman transform and master equation.  Then there is a bijective correspondence:
\begin{equation*}
Hom(\FT(\op{P}),\op{O})\cong ME(lim_\V(\op{P} \tensor\op{O}))
\end{equation*}
\end{thm}

\section{Model structures, resolutions and the W--constructions}
\label{modelpar}
In this section we discuss Quillen model structures for $\fopcat_\C$. It turns out that these model structures can be defined if $\C$ satisfies certain conditions and if this is the case work for all $\FF$, e.g.\ all the previous examples.

\subsection{Model structure}
\begin{thm}\label{modelthm}  Let $\FF$ be a Feynman category and let $\op{C}$ be a cofibrantly generated model category and a closed symmetric monoidal category having the following additional properties:
\begin{enumerate}
\item  All objects of $\op{C}$ are small.
\item  $\op{C}$ has a symmetric monoidal fibrant replacement functor.
\item  $\op{C}$ has $\tensor$-coherent path objects for fibrant objects.
\end{enumerate}
Then $\fopsc$ is a model category where a morphism $\phi\colon \op{O}\to\op{Q}$ of $\F$-\opers is a weak equivalence (resp. fibration) if and only if $\phi\colon \op{O}(v)\to\op{Q}(v)$ is a weak equivalence (resp. fibration) in $\op{C}$ for every $v\in \V$.
\end{thm}

\subsubsection{Examples}
\begin{enumerate}
\item Simplicial sets. (Straight from Theorem \ref{modelthm})
\item $dgVect_k$ for $char(k)=0$ (Straight from Theorem \ref{modelthm})
\item $Top$ (More work, see below.)
\end{enumerate}

\subsubsection{Remark}
Condition (i) is not satisfied for $Tap$ and so we can not directly apply the theorem. In \cite{feynman} this  this point was first cleared up by following \cite{Fresse} and using the fact that all objects in $Top$ are small with respect to topological inclusions.

\begin{thm}  Let $\op{C}$ be the category of topological spaces with the Quillen model structure.  The category $\fopsc$ has the structure of a cofibrantly generated model category in which the forgetful functor to $\vseq$ creates fibrations and weak equivalences.
\end{thm}

\subsection{Quillen adjunctions from morphisms of Feynman categories}
\subsubsection{Adjunction from morphisms}
We assume $\op{C}$ is a closed symmetric monoidal and model category satisfying the assumptions of Theorem \ref{modelthm}.  Let $\fr{E}$ and $\FF$ be Feynman categories and let $\alpha\colon \fr{E}\to\FF$ be a morphism between them.  This morphism induces an adjunction
\begin{equation*}
\alpha_*\colon\eops \leftrightarrows \fopsc\colon \alpha^*
\end{equation*}
where $\alpha^*(\op{A}):= \op{A}\circ\alpha$ is the right adjoint and $\alpha_*(\op{B}):=Lan_\alpha(\op{B})$ is the left adjoint.

\begin{lemma}\label{qalem} Suppose $\alpha_R$ restricted to $\vfmods\to \vemods$ preserves fibrations and acyclic fibrations, then the adjunction  $(\alpha_L, \alpha_R)$ is a Quillen adjunction.
\end{lemma}

\subsection{Example}

\begin{enumerate}
\item Recall that $\CCyclic$ and $\modular$ denote the Feynman categories whose $\opers$ are cyclic and modular operads, respectively, and that there is a morphism $i\colon\CCyclic\to\modular$ by including $\ast_S$ as genus zero $\ast_{S,0}$.
\item This morphism induces an adjunction between cyclic and modular operads
\begin{equation*}
i_*\colon\CCyclic\text{-}\op{O}ps_\op{C} \leftrightarrows  \modular\text{-}\op{O}ps_\op{C}\colon i^*
\end{equation*}
and the left adjoint is called the modular envelope of the cyclic operad.
\item  The fact that the morphism of Feynman categories is inclusion means that $i_R$ restricted to the underlying $\V$-modules is given by forgetting, and since fibrations and weak equivalences are levelwise, $i_R$ restricted to the underlying $\V$-modules will preserve fibrations and weak equivalences.
\item Thus by  the Lemma above this adjunction is a Quillen adjunction.

\end{enumerate}

\subsection{Cofibrant replacement}


\begin{thm}  The Feynman transform  of a non-negatively graded dg $\F$-$\oper$ is cofibrant.
\label{cofcor}

The double Feynman transform of a non-negatively graded dg $\F$-$\oper$ in a quadratic Feynman category is a cofibrant replacement.
\end{thm}

\subsection{W-construction}
\label{wpar}
\subsubsection{Setup}
In this section we start with a quadratic Feynman category $\FF$.

\subsubsection{The category   $w(\FF,Y)$, for  $Y \in \F$}

\mbox{}\\

{\sc Objects:}
The objects are the set $\coprod_n C_n(X,Y)\times [0,1]^n$, where $C_n(X,Y)$ are chains of morphisms from $X$ to $Y$ with  $n$ degree $\geq1$ maps
modulo contraction of isomorphisms.

  An object in $w(\FF,Y)$ will be represented (uniquely up to contraction of isomorphisms) by a diagram
\begin{equation*}
X\xrightarrow[f_1]{t_1} X_1\xrightarrow[f_2]{t_2} X_2\to\dots\to X_{n-1}\xrightarrow[f_n]{t_n} Y
\end{equation*}
where each morphism is of positive degree and where $t_1,\dots,t_n$ represents a point in $[0,1]^n$.  These numbers will be called weights.  Note that in this labeling scheme isomorphisms are always unweighted.


{\sc Morphisms:}
\begin{enumerate}
\item  Levelwise commuting isomorphisms which fix $Y$, i.e.:
\begin{equation*}
\xymatrix{X \ar[r] \ar[d]^{\cong} & X_1 \ar[d]^{\cong} \ar[r] & X_2 \ar[d]^{\cong} \ar[r] & \dots \ar[r] & X_n \ar[d]^{\cong} \ar[r] & Y  \\ X^{\prime} \ar[r] & X^{\prime}_1\ar[r] & X^{\prime}_2\ar[r] &\dots \ar[r] & X^{\prime}_n \ar[ur] & }
\end{equation*}
\item  Simultaneous $\SS_n$ action.
\item  Truncation of $0$ weights: morphisms of the form $(X_1\stackrel{0}\to X_2\to\dots\to Y)\mapsto (X_2\to\dots\to Y)$.
\item  Decomposition of identical weights:  morphisms of the form $(\dots \to X_i\stackrel{t}\to X_{i+2} \to \dots) \mapsto (\dots\to X_i\stackrel{t}\to X_{i+1}\stackrel{t}\to X_{i+2}\to\dots)$ for each (composition preserving) decomposition of a morphism of degree $\geq 2$ into two morphisms each of degree $\geq 1$.
\end{enumerate}

\begin{definition}
Let $\op{P}\in\fopst$.  For $Y \in ob(\F)$ we define
\begin{equation*}
W(\op{P})(Y):= colim_{w(\FF,Y)}\op{P}\circ s (-)
\end{equation*}
\end{definition}

\begin{thm}  Let $\FF$ be a simple Feynman category and let $\op{P}\in\fopst$ be $\rho$-cofibrant.  Then $W(\op{P})$ is a cofibrant replacement for $\op{P}$ with respect to the above model structure on $\fopst$.
\end{thm}

Here ``simple'' is a technical condition satisfied by all graph examples.

\section{Geometry}
\label{geopar}
\subsection{Moduli space geometry}
Although many of the examples up to now have been algebraic or combinatorial in nature, there are very important and deep links to the geometry of moduli spaces.
We will discuss these briefly.

\subsubsection{Modular Operads}
The typical topological example for modular operads are the Deligne--Mumford compactifications  $\bar M_{gn}$ of Riemann's moduli space of curves of genus $g$ with $n$ marked points.

These give rise to chain and homology operads. An important application comes from enumerative geometry.
 Gromov--Witten invariants make $H^*(V)$ an algebra over $H_*(\bar M_{g,n})$ \cite{ManinBook}.

\subsubsection{Odd Modular}
As explained in \cite{KWZ}, the canonical geometry for odd modular operads is given by $\bar M^{KSV}$ which are real blowups of $\bar M_{gn}$ along the boundary divisors \cite{KSV}.

On the topological level one has $1$-parameter gluings parameterized by $S^1$. Taking the full $S^1$ family on chains or homology gives
us the structure of an odd modular operad. That is the gluing operations have degree $1$ and in the dual graph, the edges have degree $1$.

\subsection{Master Equation and compactifications}
 Going back to Sen and Zwiebach \cite{Zwie}, a viable string field theory action $S$ is a solution of the quantum master equation.
 Rephrasing this one can say ``The master equation drives the compactification'', which is one of the mantras of \cite{KWZ}.

In particular, the constructions of \cite{KSV} and \cite{HVZ} give the correct compactification.

\subsection{W--construction}
In \cite{BK} we will prove the fact that the derived modular envelope defined via the $W$--construction of the cyclic associative operads
is the Kontsevich/Penner compactification $M^{comb}_{g,n}$.

We will also give an $A_{\infty}$ version of this theorem and a 2--categorical realization that gives our construction of string topology and  Hochschild operations from Moduli Spaces \cite{hoch1,hoch2}
via the Feynman transform.

\section{Bi-- and Hopf algebras}
\label{Hopfpar}
We will give a brief overview of the constructions of \cite{GKT}.

\subsection{Overview}
Consider  a non--Sigma Feynman category $\B=Hom(Mor(\F),\Z)$ .

{\sc Product.}
Assume that $\FF$ is strict monoidal, that is $\F$ is strict monoidal, then $\otimes$ is an associative unital product on $\B$ with unit $id_{\Eins_\F}$.

{\sc Coproduct.} Assume that $\F$ decomposition finite, i.e.\ that the sum below is finite.
Set
\begin{equation}\Delta(\phi)=\sum_{(\phi_0,\phi_1):\phi=\phi_1\circ \phi_0 }\phi_0\otimes \phi_1\end{equation}
and $\eps(\phi)=1$ if $\phi=id_X$ and $0$ else.

\begin{thm}\cite{GKT}
$\B$ together with the structures above is a bi--algebra. Under certain mild assumptions, a canonical quotient is a Hopf algebra.
\end{thm}

\begin{rmk} Now, it is not true that any strict monoidal category with finite decomposition yields a bi--algebra. Also, if $\FF$ is a Feynman category,
then $\FF^{op}$, although not necessarily a Feynman category, does yield a bi--algebra.
\end{rmk}

\subsubsection{Examples}
The Hopf algebras of Goncharov for multi--zeta values \cite{Gont} can be obtained in this way starting with the Joyal dual of the surjections in the augmented simplicial category. In short, this Hopf algebra structures follows from the fact that simplices form an operad.
 In a similar fashion, but using a graded version, we recover a Hopf algebra of Baues that he defined for double loop spaces \cite{Baues}.
 We can also recover the non--commutative Connes--Kreimer Hopf algebra of planar rooted trees, see e.g.\ \cite{foissyCR1} in this way.

\begin{rmk}
This coproduct for any 
finite decomposition category appeared in
\cite{Moebiusguy} and was picked up later in \cite{JR}.
We realized with hindsight that the co--product we first constructed on indecomposables,  as suggested to us by Dirk Kreimer, is equivalent to this coproduct.
\end{rmk}

\subsubsection{Symmetric version}

There is a version for symmetric Feynman categories, but the constructions are more involved.
In this fashion, we can reproduce Connes--Kreimer's Hopf algebra. There is a three-fold hierarchy. A bialgebra version, a commutative Hopf algebra version and an ``amputated" version, which is actually the algebra considered in \cite{CK}. A similar story holds for the graph versions and in general.

\subsection{Details: Non--commutative version}
We use non--symmetric  Feynman categories whose underlying tensor structure is only monoidal (not symmetric). $\V^{\otimes}$ is the free monoidal category.

\begin{lem}[Key Lemma]
The bi--algebra equation holds due to the hereditary condition (ii).
\end{lem}
The proof is a careful check of the diagrams that appear in the bialgebra equation.

For  $\Delta\circ\mu$ the sum is over diagrams of the type
\begin{equation}
\label{deltamueq}
\xymatrix{
X\otimes X'\ar[rr]^{\Phi=\phi\otimes \psi}\ar[dr]^{\Phi_0}&&Z\otimes Z'\\
&Y\ar[ur]^{\Phi_1}&\\
}
\end{equation}
where $\Phi=\Phi_1\circ\Phi_0$.

When considering $(\mu\otimes \mu)\circ\pi_{23}\circ (\Delta\otimes \Delta)$
the diagrams are of the type
\begin{equation}
\label{mudeltaeq}
\xymatrix{
X\otimes X'\ar[rr]^{\phi\otimes \psi}\ar[dr]^{\phi_0\otimes \psi_0}&&Z\otimes Z'\\
&Y\otimes Y'\ar[ur]^{\phi_1\otimes \psi_1}&\\
}
\end{equation}
where $\phi=\phi_1\circ\phi_0$
and
$\psi=\psi_1\circ\psi_0$. In general, there is no reason for there to be a bijection of such diagrams, but there is for non--symmetric Feynman categories.

For simplicity, we assume that $\F$ is  skeletal.

\subsection{Hopf quotient}
Even after quotienting out by the isomorphisms, the bi--algebra is usually not connected. The main obstruction is that there are many identities and that there are still automorphisms. The main point is that in the skeletal case:
\begin{equation}
\Delta(id_X)=\sum_{\sigma \in Aut(X)}\sigma\otimes \sigma^{-1}
\end{equation}
where here and in the following we assume that if $\sigma$ has a one--sided inverse then it is invertible. This is the case in all examples.

\subsubsection{Almost connected Feynman categories}
In the skeletal version, consider the ideal generated by ${\mathfrak C}=|Aut(X)|[id_X]-|Aut(Y)|[id_Y]\subset B$, this is closed under $\Delta$, but not quite a co--ideal.
Rescaling $\eps$ by $\frac{1}{|Aut(X)|}$, $\H=\B/{\mathfrak C}$ becomes a bi--algebra.
We  call $\FF$ almost connected if $\H$ is connected.

\begin{thm}
For the almost connected version $\H$ is
a connected bi--algebra and hence a Hopf--algebra.
\end{thm}

\subsection{Symmetric/Commutative version}
In the case of a symmetric Feynman category, the bi--algebra equation does not hold anymore, due to the fact that $Aut(X)\otimes Aut(Y)\subset Aut(X\otimes Y)$
may be a proper subgroup due to the commutativity constraints. The typical example is $\alephsym$ where
$Aut(n)\times Aut(m)=\SS_n\times \SS_m$ while $Aut(n+m)=\SS_{n+m}$.
In order to rectify this, one considers the co--invariants. Since commutativity constraints are
isomorphisms the resulting algebra structure is commutative.

Let $\B^{iso}$ the quotient by the  ideal  defined by the equivalence relation generated by isomorphism. That is $f\sim g$ if there are isomorphisms $\sigma,\sigma'$ such that $f=\sigma\circ g\circ \sigma'$. This ideal is again closed under co--product.
As above one can modify the co--unit to obtain a bialgebra structure on $\B^{iso}$. Now the ideal generated by ${\mathfrak C}=\langle |Aut(X)|[id_X]-|Aut(Y)|[id_Y]$ is a co--ideal and  $\H=\B/{\mathfrak C}$ becomes a bi--algebra. We  call $\FF$ almost connected if $\H$ is connected.

The main theorem is
\begin{thm}
If $\FF$ is almost connected, the coinvariants $\B^{iso}$ are a commutative Hopf algebra.
\end{thm}

This allows one to construct Hopf algebras with external legs in the graph examples. It also explains why the Connes--Kreimer examples are commutative.

\subsubsection{Amputated version}  In order to forget the leg structure, aka.\ amputation, one needs a semi--cosimplicial structure, i.e. one must be able to forget external legs coherently. This is always possible by deleting flags in the graph cases.
Then there is a colimit, in which all the external legs can be forgotten.
Again, one obtains a Hopf algebra. The example {\em par excellence} is of course,
Connes--Kreimer's Hopf algebra without external legs (e.g. the original  version).

\subsection{Restriction and Generalization of special case: co--operad with multiplication}
In a sense the above examples were free. One can look at a more general setting where this is not the case.
This is possible in the simple cases of enriched Feynman categories over $\Surj$. Here the morphisms are operads, and $\B$ has the dual co--operad structure for the one--comma generators. The tensor product $\otimes$ makes $\B$ have the structure of a free algebra over the one--comma generators $\O(n)$ with the co--operad structure being distributive or multiplicative over $\otimes$.
Now one can generalize to a general co--operad structure with multiplication.

\subsubsection{Coproduct for a cooperad with multiplication}

\begin{thm}\cite{GKT}
\label{bialgthm}
Let $\Coop$ be a co-operad with compatible associative multiplication. $\mu:\Coop(n)\otimes\Coop(m)\to \Coop(n+m)$ in an Abelian symmetric monoidal category with unit $\unit$. Then
$
\B:=
\bigoplus_n \Coop(n)
$
is a (non-unital, non-co-unital) bialgebra, with multiplication $\mu$ and comultiplication $\Delta$ given by $(\id\otimes\mu)\check\gamma$:
{\footnotesize
\begin{equation}
\label{deltadefeq}
\vcenter{\xymatrix@R-1.3pc{\ar@/_6mm/[rrdd]_(0.3){\textstyle{\Delta:=(\id\otimes\mu)\check\gamma}\;\;}\Coop(n)
\rrto^-{\;\check\gamma\;}
&&\hspace*{-6mm}
{\displaystyle\bigoplus_{\substack{k\geq1,\\ n=m_1+\dots+m_k}}}
\hspace*{-6mm}
\left(\Coop(k)\otimes{\displaystyle\bigotimes_{r=1}^k\Coop(m_r)}\right)
\ar[]-<0mm,6mm>;[dd]^(0.4){\;\id{}\,\otimes\,{}\mu^{k-1}\;}\\ \\
&&\displaystyle\bigoplus_{k\geq1}
\Coop(k)\otimes\Coop(n).
}}
\end{equation}
}
\end{thm}

\subsubsection{Free cooperad with multiplication on a cooperad} The guiding example is:
$$\Coop^{nc}(n)=\bigoplus_k\bigoplus_{(n_1,\dots,n_k):\sum n_i=n} \Coop(n_1)\odo \Coop(n_k)$$
Multiplication is given by $\mu=\otimes$. This structure coincides with  one of the constructions of a non--connected operad in \cite{KWZ}.

The example is the one that is relevant for the three Hopf algebras of Baues, Goncharov and Connes--Kreimer.
It also shows how a cooperad with multiplications generalizes  an enrichment of $F_{surj}$.

This is most apparent in Connes--Kreimer, where the Hopf algebra is not actually on rooted trees, but rather on forests.
The extension of the co--product to a forest is tacitly given by the bi--algebra equations.

In the symmetric case, one has to further induce the natural $(\SS_{n_1}\times \dots \times \SS_{n_k})\wr \SS_k$ action to an $\SS_n$ action for each summand. The coinvariants constituting $\B^{iso}$
are then the symmetric products $ \Coop(n_1)_{\SS_{n_1}}\odot \dots \odot  \Coop(n_k)_{\S_{n_k}}$.

The following is the list of motivating examples:

\mbox{}

\begin{tabular}{l|l|ll}
Hopf algebras&(co)operads&Feynman category\\
\hline
$H_{Gont}$&$Inj_{*,*}=Surj^*$&$\FF_{Surj}$\\
$H_{CK}$&leaf labelled trees&$\FF_{Surj,\O}$\\
$H_{CK,graphs}$&graphs&$\FF_{graphs}$\\
$H_{Baues}$&$Inj_{*,*}^{gr}$&$\FF_{Surj,odd}$
\end{tabular}

\subsubsection{Grading/Filtration, the $q$ deformation and infinitesimal version}
We will only make very short remarks, the details are in \cite{GKT}.

The length of an object in the Feynman category setting is  replaced  by a depth filtration. The algebras are then deformations of their associated graded, see \cite{GKT}.
In the amputated version one has to be more careful with the grading.

\mbox{}

\begin{tabular}{l|l}
Co-operad with multiplication&operad degree $-$ depth\\
Amputated version&co-radical degree $+$ depth\\
\end{tabular}

\mbox{}

Taking a slightly different quotient, one can get a non--unital, co--unital bi--algebra and a $q$--filtration. Sending $q\to 1 $ recovers $\mathcal H$.

%
%
%

\bibliography{lectbib}
\bibliographystyle{halpha}

\end{document}